\documentclass[12pt]{article}
\usepackage{amsthm,amsfonts,amssymb,amscd}

\begin{document}

\begin{center}
\Large \bf Birationally rigid iterated Fano double covers
\end{center}
\vspace{1cm}

\centerline{\large  \bf A.V.Pukhlikov\footnote{The author was financially supported by Russian 
Foundation of Basic Research, grants 02-15-99258 and
02-01-00441, by a Science Support Foundation grant for young researchers and by INTAS-OPEN, grant 2000-269.}}\vspace{1cm}

\begin{center}
Max-Planck-Institut f\" ur Mathematik \\
Vivatsgasse 7 \\
53111 Bonn \\
GERMANY \\
e-mail: {\it pukh@mpim-bonn.mpg.de}
\end{center}
\vspace{0.5cm}

\begin{center}
Steklov Institute of Mathematics \\
Gubkina 8 \\
117966 Moscow \\
RUSSIA \\
e-mail: {\it pukh@mi.ras.ru}
\end{center}
\vspace{1cm}

\centerline{November 17, 2002}
\vspace{1cm}

\parshape=1
3cm 10cm \noindent {\small \quad\quad\quad \quad\quad\quad\quad
\quad\quad\quad {\bf Abstract}\newline
Iterating the procedure of making a double cover over a given variety, we construct large families of smooth higher-dimensional Fano varieties of index 1. These varieties can be realized as complete intersections in various weighted projective spaces. A generic variety in these families is proved to be birationally superrigid; in particular, it admits no non-trivial structures of a fibration into rationally connected (or uniruled) varieties, it is non-rational and its groups of birational and biregular self-maps coincide.} 

\newpage

{\small

CONTENTS
\vspace{0.7cm}

\noindent 0. Introduction

0.1. Birationally rigid varieties

0.2. The known results and natural conjectures

0.3. The iterated Fano double covers

0.4. The structure of the paper 

0.5. Historical remarks and acknowledgements
\vspace{0.5cm}

\noindent 1. The method of maximal singularities

1.1. A criterion of birational superrigidity

1.2. The first proof: counting multiplicities

1.3. The second proof: the connectedness principle

1.4. The Lefschetz theorem

1.5. Proof of Proposition 0.1
\vspace{0.5cm}

\noindent 2. Iterated double covers

2.1. Coordinate presentations

2.2. The regularity condition

2.3. Non-regular sets of polynomials

2.4. Start of the proof of Proposition 2.1

2.5. The polynomials $h_{i,j}$ depend on each other

2.6. Estimates for the codimension
\vspace{0.5cm}

\noindent 3. Hypertangent divisors

3.1. How to obtain a bound for the multiplicity

3.2. Construction of hypertangent divisors

3.3. The regularity condition for hypertangent divisors

3.4. The Lefschetz theorem once again

3.5. The points of class $e=0$

3.6. The points of class $e\geq 1$
\vspace{0.5cm}

\noindent References }

\newpage

\section*{Introduction}

\subsection{Birationally rigid varieties}

Birational rigidity is one of the most striking phenomena of higher-dimensional algebraic geometry. Speaking informally,  birational rigidity means that certain algebraic varieties, on which there are no non-zero global regular differential forms (that is, rationally connected varieties), behave as if they were of general type. Starting from the pioneer paper of V.A.Iskovskikh and Yu.I.Manin on the three-dimensional quartic [IM] of 1970 it has been gradually understood that in higher dimensions birationally rigid varieties and fibrations not only do not form an exceptional exotic class,  but on the contrary, are quite typical. Sarkisov's theorem [S1,S2] means that in a sense ``almost all'' conic bundles of dimension three and higher are birationally rigid. In [P1,P7,P8,P10,P12] it was proved that Fano hypersurfaces and more generally certain complete intersections of index 1 are birationally rigid. In [P6,P9] it was proved that ``almost all'' del Pezzo fibrations and fibrations into Fano hypersurfaces over ${\mathbb P}^1$ are birationally rigid. Nowadays it is clear that, on the contrary, non-rigid Fano varieties and fibrations are less typical.

In the present paper we use the more traditional definition of birational rigidity. We work over the field ${\mathbb C}$ of complex numbers. A smooth variety $X$ of dimension
$\geq 3$ is said to be
{\it birationally superrigid}, if for each birational map
$\chi\colon X-\,-\,\to X'$ onto a variety $X'$ of the same
dimension, smooth in codimension one, and each linear system
$\Sigma'$ on $X'$, free in codimension 1 (that is,
$\mathop{\rm codim}\nolimits \mathop{\rm Bs}\Sigma'\geq 2$),
the inequality
\begin{equation}
\label{i1}
c(\Sigma,X)\leq c(\Sigma',X')
\end{equation}
holds, where $\Sigma=(\chi^{-1})_*\Sigma'$ is the proper inverse image
of $\Sigma'$ on $X$ with respect to $\chi$, and $c(\Sigma,X)=c(D,X)$
stands for the {\it threshold of canonical adjunction}
$$c(D,X)=\sup\{b/a|b,a\in{\mathbb Z}_+\setminus \{0\},
|aD+bK_X|\neq\emptyset\}$$
$D\in\Sigma$, and similarly for $\Sigma'$, $X'$. \par

A smooth variety $X$ of dimension $\geq 3$ is said to be
{\it birationally rigid}, if for each $X'$, $\chi$ and
$\Sigma'$ as above there exists a birational {\it self}-map
$$
\chi^*=\chi^*_{(X',\chi,\Sigma')}\in
\mathop{\rm Bir} X,
$$
depending on the triple $(X',\chi,\Sigma')$, such that
$$
c(\Sigma^*,X)\leq c(\Sigma',X'),
$$
where $\Sigma^*$ is the strict transform of $\Sigma$
with respect to $(\chi^*)^{-1}$, or, equivalently, the
strict transform of $\Sigma$ with respect to 
the composition 
$$
\chi\circ\chi^* \colon X -\,-\,\to X'.
$$

The difference between the rigid and superrigid cases
is not too big. Roughly speaking, superrigidity is
rigidity combined with the property that the groups of
birational and biregular maps coincide (indeed, if this
is the case, then twisting by a birational$=$biregular
self-map $\chi^*$ does not change the threshold of
canonical adjunction, hence rigidity implies the
inequality (\ref{i1}), that is, superrigidity).

The most natural object for the rigidity theory is
formed by Fano varieties. Assume that $X$ is a smooth
Fano variety of dimension $\geq 3$ with 
$\mathop{\rm Pic} X={\mathbb Z}K_X$. The important
geometric properties of birationally rigid and
superrigid varieties from this class (the properties
that motivate the very choice of the word ``rigidity'')
are summarized in the following
\vspace{0.3cm}

{\bf Proposition 0.1.} {\it Assume that $X$ is rigid. Then:
\vspace{0.1cm}

{\rm (i)} $X$ can not be fibered into uniruled varieties by a
non-trivial rational map,
\vspace{0.1cm}

{\rm (ii)} if $\chi\colon X -\, -\, \to X'$ is a birational map
onto a Fano variety $X'$  with ${\mathbb Q}$-factorial terminal
singularities such that
$\mathop{\rm Pic} X'\otimes{\mathbb Q}={\mathbb Q}K_{X'}$, then
$X$ is (biregularly) isomorphic to $X'$ (although $\chi$
itself is not necessarily an isomorphism),
\vspace{0.1cm}

{\rm (iii)} $X$ is non-rational.
\vspace{0.1cm}

\noindent Assume moreover that $X$ is superrigid. Then any
 birational map
onto a Fano variety $X'$  with ${\mathbb Q}$-factorial terminal
singularities such that
$\mathop{\rm Pic} X'\otimes{\mathbb Q}={\mathbb Q}K_{X'}$
is a (biregular) isomorphism. In particular, the
groups of birational and biregular self-maps coincide:
$$
\mathop{\rm Bir}X=\mathop{\rm Aut}X.
$$
}
\vspace{0.3cm}

These implications of
birational rigidity are well known. For convenience of the reader, we give a
(very easy) proof of Proposition 0.1 below in Sec. 1.5.

Note that Corti and Reid [C2,CR] take the properties (i) and (ii)
(property (iii) is an immediate implication of (i)) as a
{\it definition} of birational rigidity.

Although in this paper we study only Fano varieties, to make the picture
complete, let us give the definition of birational rigidity for Fano
fibrations, too. The relative version is completely analogous to the
absolute one (see the definition above), with the only difference: the 
group of birational self-maps is replaced by the group of {\it fiber-wise}
birational self-maps.

Assume that $X/S$ is a rationally connected fibration:
$$
\begin{array}{rcl}\displaystyle
   &   X   &   \\  \displaystyle
\pi  &  \downarrow & \mbox{generic fiber}\,\, F_{\eta} \\
\displaystyle    &    S   &   
\end{array}
$$
The fiber of general position is assumed to be rationally
connected, so that $X$ itself is uniruled and the thresholds
of canonical adjunction are finite. We define the group of
{\it proper} birational self-maps of the fibration $X/S$,
setting
$$
\mathop{\rm Bir} (X/S)=\mathop{\rm Bir} F_{\eta},
$$
that is, $\mathop{\rm Bir} (X/S)\subset\mathop{\rm Bir} X$
is the subgroup consisting of all maps 
$\chi\colon X-\,-\,\to X$
such that $\chi$ transforms each fiber into itself. Now
the fibration is said to be {\it birationally rigid}
(as a fibration!), if for any variety $X'$ of the same
dimension, smooth in codimension one, any birational map  $\chi\colon X-\,-\,\to X'$ and any linear system $\Sigma'$
on $X'$ there exists a self-map
$$
\chi^*\in\mathop{\rm Bir} (X/S)
$$
such that for the strict transform $\Sigma^*$ of the
linear system $\Sigma'$ with respect to the composition
$(\chi\circ\chi^*)^{-1}$ the inequality
$$
c(\Sigma^*,X)\leq c(\Sigma',X')
$$
is satisfied.
\vspace{0.3cm}

{\bf Proposition 0.2.} {\it Assume that $X/S$ is a Fano
fibration with $X$, $S$ smooth such that
$$
\mathop{\rm Pic} X={\mathbb Z}K_X\oplus
\pi^* \mathop{\rm Pic} S
$$
and for any effective class $D=mK_X+\pi^*T$ the
class $NT$ is effective on $S$ for some
$N\geq 1$. Assume furthermore
that $X/S$ is birationally rigid. Then for any
rationally connected fibration $X'/S'$ and any
birational map
$$
\chi\colon X -\,-\,\to X'
$$
(provided that such maps exist) there is a rational
dominant map
$$
\alpha\colon S -\,-\,\to S'
$$
making the following diagram commutative:
$$
\begin{array}{rcccl}
    &  X  &  \stackrel{\chi}{-\,-\,\to} & X' &  \\
\pi  &  \downarrow &    &  \downarrow  & \pi'  \\
    &  S   & \stackrel{\alpha}{-\,-\,\to} &  S'.  &
\end{array}
$$
}


\subsection{The known results and natural conjectures}

At the moment birational superrigidity is proved for the following
classes of higher-dimensional Fano varieties:

\begin{itemize}

\item
smooth hypersurfaces $V_M\subset{\mathbb P}^M$ of degree $M$, $M\geq 4$
(for $M=4$ it is the classical theorem of V.A.Iskovskikh and Yu.I.Manin [IM],
the case $M=5$ was proved in [P1], for generic hypersurfaces of degree
$M\geq 6$ superrigidity was proved in [P7], finally, for arbitrary smooth
hypersurfaces it was proved in [P12]);

\item
generic (in the sense of Zariski topology) Fano complete intersections
$$ 
\begin{array}{rcl} \displaystyle
V= & V_{d_1\cdot d_2\cdot \dots \cdot d_k} &
\subset {\mathbb P}^{M+k}\\ \displaystyle
& \parallel & \\ \displaystyle &
F_1\cap F_2\cap \dots \cap F_k, &
\end{array}
$$
where $\mathop{\rm deg} F_i=d_i$, $d_1+\dots+d_k=M+k$, $d_i\geq 2$
and the inequality $2k<M$ holds [P10];

\item
generic (in the sense of Zariski topology) Fano double covers
$$
\sigma\colon V\to Q=Q_m\subset {\mathbb P}^{M+1},
$$
where $M\geq 4$, $Q$ is a hypersurface of degree $m\geq 3$, $\sigma$ is
ramified over the divisor
$$
W=W_{m\cdot 2l}=W^*\cap Q,
$$
where $W^*=W^*_{2l}\subset {\mathbb P}^{M+1}$ is a hypersurface of
degree $2l$, $m+l=M+1$ [P8]. For $m=1,2$ birational superrigidity was earlier
proved in [P2] for smooth double spaces and double quadrics without the assumption of them being generic;

\item
the known examples [P3,P11] show that mild singularities do not change the
picture: the varieties remain rigid, and usually even superrigid.

\end{itemize}
\vspace{0.1cm}

{\bf Conjecture 1 (the absolute case).} {\it A smooth
Fano variety $V$ of dimension $\mathop{\rm dim}V\geq 5$
with the Picard group $\mathop{\rm Pic}V={\mathbb Z} K_V$
is birationally superrigid.}
\vspace{0.3cm}

The assumption of $V$ being smooth seems to be unnecessarily strong. 
\vspace{0.3cm}

{\bf Conjecture 2 (the absolute case).} {\it A Fano variety $V$ of dimension $\mathop{\rm dim}V\geq 4$ with factorial terminal singularities and the Picard group $\mathop{\rm Pic}V={\mathbb Z} K_V$ is birationally rigid.}
\vspace{0.3cm}

One can make Conjecture 2 stronger in its turn, replacing the factorial singularities by the ${\mathbb Q}$-factorial ones and modifying the condition on the Picard group in an appropriate way. However, Fano varieties with terminal quotient singularities have not yet been studied in higher ($\geq 4$) dimensions from the viewpoint of their birational rigidity. Since there are no completely studied examples it seems it is yet too early to formulate general conjectures.

To complete the picture, we give a relative version of the conjectures on birational rigidity.

Let $V\to S$ be a fibration into Fano varieties satisfying the following standard conditions:

(i) $V$ is smooth, $\mathop{\rm dim} V\geq 4$,

(ii) $\mathop{\rm dim} S\geq 1$, the anticanonical class $-K_V$ is relatively ample,

(iii) $\mathop{\rm Pic} V={\mathbb Z} K_V\oplus \pi^*\mathop{\rm Pic} S$.

(In dimension $\mathop{\rm dim} V=3$ this means that $V/S$ is a Mori fiber space.)
\vspace{0.3cm}

{\bf Conjecture 3 (the relative case).} {\it  If the fibration $V/S$ is sufficiently twisted over the base $S$, then it is birationally superrigid.}
\vspace{0.3cm}

If $\mathop{\rm dim} S=1$, that is, $S={\mathbb P}^1$, then the twistedness of the fibration $V/S$ is characterized by the properties of the numerical class
$$
K^2_V\in A^2(V),
$$
more precisely, in the known cases it is sufficient to assume
that the following ``$K^2_V$-condition'' is satisfied: $K^2_V$ does not lie in the interior of the cone of the effective cycles of codimension two [P6,P9]. However, this $K^2_V$-condition can be somewhat weakened, see [G1,G2,Sob1,Sob2].

In the general case the ``twistedness'' of the fibration $V/S$ can be imagined in the following way. Let the fiber of general position $F_s$, $s\in S$, belong to a fixed family of Fano varieties ${\cal F}$. For simplicity we assume that for a general member $F\in{\cal F}$ the anticanonical class $-K_V$ is very ample and determines an embedding
$$
F\hookrightarrow {\mathbb P}^N.
$$
Let ${\cal H}=\mathop{\rm Hilb}(F)$ be the Hilbert scheme of {\it embedded} varieties of the family ${\cal F}$. Then the fibrations $V/S$ with the general fiber in ${\cal F}$ correspond to the maps
$$
S\to{\cal H}.
$$
The degree of the image of $S$ with respect to a fixed class on ${\cal H}$ can be considered as a characteristic of twistedness of $V/S$. However, for particular fibrations the twistedness condition takes a simple form and is easy to check, see [P6,P9].


\subsection{Iterated Fano double covers}

The aim of the present paper is to prove birational superrigidity of a large class of higher-dimensional Fano varieties, generalizing Fano double hypersurfaces [P8].
The varieties under consideration can be realized as complete
intersections in weighted projective spaces. Let us give their construction.

The ground field is assumed to be the field of complex numbers ${\mathbb C}$. The symbol ${\mathbb P}$ stands for the projective space ${\mathbb P}^{M+k}$, where $M\geq 4$, $2k\leq M-1$. 
Let us choose a system of homogeneous coordinates on ${\mathbb P}$,
say $(x_0:x_1:\dots:x_{M+k})$. Consider a set of homogeneous
polynomials
$$
f_1,\dots,f_k,g_1,\dots,g_m
$$
in the variables $(x_*)$ of degrees, respectively,
$$
d_1,\dots,d_k,2l_1,\dots,2l_m,
$$
where $m\geq 1$, $l_i\geq 2$ and the following equality
holds:
$$
\sum^k_{i=1}d_i+\sum^m_{i=1}l_i=M+k.
$$
Set
$$
Q(f_*)=Q(f_1,\dots,f_k)=
\mathop{\bigcap}\limits^k_{i=1}\{f_i=0\}
$$
to be the complete intersection of hypersurfaces
$F_i=\{f_i=0\}$ in ${\mathbb P}$. Set also
$$
W_i=\{g_i=0\}\subset{\mathbb P}.
$$
Define a sequence of double covers
$$
\sigma_i\colon {\mathbb P}^{(i)} \to {\mathbb P}^{(i-1)},
$$
where ${\mathbb P}^{(0)}={\mathbb P}$, $i=1,\dots,m$, in the following way. 
The cover $\sigma_1$ is branched over $W_1$. Assume that $\sigma_1,\dots,\sigma_i$ are already constructed. For an
arbitrary $j\in\{0,\dots,i-1\}$ set
$$
\sigma_{i,j}=\sigma_{j+1}\circ\dots\circ\sigma_i\colon
{\mathbb P}^{(i)}\to{\mathbb P}^{(j)}.
$$
In particular, $\sigma_i=\sigma_{i,i-1}$. Obviously $\sigma_{i,j}$ is a finite morphism of degree $2^{i-j}$.
Now the double cover
$$
\sigma_{i+1}\colon {\mathbb P}^{(i+1)}\to{\mathbb P}^{(i)}
$$
is determined by the branch divisor
$$
{\widetilde W}_{i+1}=\sigma^{-1}_{i,0}(W_{i+1})
\subset {\mathbb P}^{(i)}.
$$
As a result we get a sequence of double covers
$$
{\mathbb P}^{(m)}\stackrel{\sigma_m}{\to}{\mathbb P}^{(m-1)}
\stackrel{\sigma_{m-1}}{\to}\dots\stackrel{\sigma_1}{\to}
{\mathbb P}^{(0)}={\mathbb P}.
$$
Finally, set
$$
\sigma=\sigma_{m,0}\colon {\mathbb P}^{(m)}\to {\mathbb P}.
$$
It is a finite morphism of degree $2^m$.

The collections of homogeneous polynomials $(f_*;g_*)$ are
parametrized by the space
$$
{\cal H}=
\prod^k_{i=1}\left[
H^0({\mathbb P}, {\cal O}_{{\mathbb P}}(d_i))\setminus \{0\})
\right] \times
\prod^m_{i=1}\left[
H^0({\mathbb P}, {\cal O}_{{\mathbb P}}(2l_i))\setminus \{0\})
\right].
$$
For a general collection $(f_*;g_*)\in{\cal H}$ all the varieties
$$
Q=Q(f_*), \quad {\mathbb P}^{(i)}, \quad {\widetilde W}_i, \quad {\widetilde W}_i\cap \sigma^{-1}_{i-1,0}(Q(f_*)),
$$
are obviously smooth. Set
$$
V=V(f_*;g_*)=\sigma^{-1}(Q)\subset {\mathbb P}^{(m)}.
$$
This is a smooth subvariety of codimension $k$ in ${\mathbb P}^{(m)}$. Obviously
$$
V=\mathop{\bigcap}\limits^{k}_{i=1}\sigma^{-1}(F_i)
$$
is a smooth complete intersection in ${\mathbb P}^{(m)}$. It is easy to see that
$$
\begin{array}{ccc} \displaystyle
K_V  & = & \displaystyle
\left[
-(M+k+1)+\sum^k_{i=1}d_i+\sum^m_{i=1}l_i
\right] \sigma^* H \\  \\  \displaystyle
   &   &  \| \\  \\ \displaystyle
   &   &   -\sigma^* H,
\end{array}
$$
where $H\in\mathop{\rm Pic} {\mathbb P}$ is the class of a hyperplane. By the Lefschetz theorem we get
$\mathop{\rm Pic} V={\mathbb Z} \sigma^* H$, so that $V$ is
a smooth Fano variety of index 1 and dimension $M$.

Here is the main result of the present paper.
\vspace{0.5cm}

{\bf Theorem.} {\it There is a non-empty open Zariski subset
$U\subset {\cal H}$ such that for any collection of polynomials
$(f_*;g_*)\in U$ the variety $V=V(f_*;g_*)$ is a (smooth)
birationally superrigid variety.}
\vspace{0.5cm}

The open subset $U$ is defined below by explicit {\it regularity conditions}. There is no doubt that {\it any smooth} variety $V=V(f_*;g_*)$ is birationally superrigid.
However, both the classical technique of hypertangent divisors, which is used in this paper, and the new technique 
based on the connectedness principle of Shokurov and 
Koll\' ar, are not strong enough to prove this fact.


\subsection{The structure of the paper}

The paper is organized in the following way. Section 1 contains the ``general theory'' of the method of maximal
singularities, which for convenience of the reader is given
here with the necessary details. In Sec. 1.1 we give a criterion of birational superrigidity, reducing the proof of
birational superrigidity to checking certain explicit conditions for subvarieties of codimension two. In Sec. 1.2 
we give the first proof of the main claim of this criterion.
This proof makes use of the technique of counting multiplicities. Sec. 1.3 contains another proof, which makes
use of Corti's idea [C2] of reduction to a non-log-canonical
singularity of a linear system on a smooth surface. The corresponding two-dimensional fact is proved by means of an
elementary technique. Actually, it is a simple implication
of one fact which was proved already in [P3].  In Sec. 1.4 the
proof of the basic criterion of Sec. 1.1 is completed.
Finally, in Sec. 1.5 we give an elementary proof of Proposition 0.1, describing geometry of birationally rigid varieties.

Section 2 is of more technical character. Its aim is to describe the regularity condition that defines the open set
$U\subset{\cal H}$ in the main theorem and to prove that this set
is non-empty. In order to do it, we use the method developed
in [P10]. However, in the present case this method should be
modified. In this section we also give certain convenient
coordinate presentations, which are used later.

Section 3 contains the heart of the proof --- there we check
the criterion of Section 1.1 for the regular iterated double
covers. In Sec. 3.1 the technique of estimating the multiplicity of a subvariety at a given point is developed.
It is
based on the concept of a hypertangent divisor. The method
is presented in a general form. After that in Sec. 3.2 we 
construct a family of hypertangent divisors for an arbitrary
point $x\in V$ of the iterated double cover $V$. Sec. 3.3
studies some properties of hypertangent divisors. These
properties follow from the regularity condition. It is here
that this condition is made use of. Finally, in Sec. 3.4 --
3.5 we complete the proof of the main theorem.


\subsection{Historical remarks and acknowledgements}

Not pretending to be complete in any sense, we just point
out the principal landmarks in the development of the theory
of birational rigidity. The very phenomenon of rigidity was
guessed by Fano [F1-F3] when he tried to extend Noether's
method [N] to dimension three. However, the techniques of
algebraic geometry of his time were not strong enough (intersection theory, sheaves and cohomology, resolution of
singularities either were non-existent at all or just
made its first steps) to enable him to obtain complete
results. The first outlines of the theory of birational
rigidity can be seen in the papers of Yu.I.Manin of 60ies
on surfaces over non-closed fields, see, for instance, [M1-M3],
where, using B.Segre's earlier results, it is proved what is now
formulated as birational superrigidity of del Pezzo surfaces
of degree 1 with the Picard group ${\mathbb Z}$. In [M2,M3] a graph is associated to a finite sequence of blow ups. Its combinatorial invariants are very important for the classical
technique of counting multiplicities.

The decisive step was made by V.A.Iskovskikh and Yu.I.Manin
in [IM] where birational superrigidity of a smooth three-dimensional quartic was established (the authors prove
the coincidence of birational and biregular isomorphisms of
three-dimensional quartics, however the arguments of the
paper do not need any modification to produce birational
superrigidity, so that essentially it is superrigidity that
was proved there). After that the technique that was
developed in [IM] was applied in [I1] to a few families of
Fano 3-folds, which resulted in proving their birational
rigidity.

The next step was made by V.G.Sarkisov in [S1,S2]. Starting 
with V.A.Iskovskikh's results on the surfaces with a pencil
of rational curves, V.G.Sarkisov proved birational rigidity
of conic bundles with a ``sufficiently big'' discriminant
divisor. The concept of ``birational rigidity'' was in 1980
yet non-existent, so in [S1,S2] the main result is formulated
as uniqueness of the conic bundle structure on a given 
variety. The breakthrough that was made in [S1,S2] was the more
impressing that embraced conic bundles in arbitrary dimension. 

In a few years after Sarkisov's  theorem the first attempts
were made to extend the three-dimensional technique of
V.A.Iskovskikh and Yu.I.Manin to the field of higher-dimensional Fano varieties [P1,P2]. Besides, birational
geometry of a three-dimensional quartic with a non-degenerate
double point was described [P3]. However, the technical side
of this work was getting more and more complicated. The
methods needed to be improved.

In [P4] and especially [P5] (the latter paper was written in
1995 and has been distributed among the experts since that
time, although published only in 2000) the classical technique of the method of maximal singularities was
essentially simplified and clarified, which in particular
made it possible to prove birational rigidity of del Pezzo
fibrations over ${\mathbb P}^1$ [P6] --- that is, the only class of three-dimensional Mori fiber spaces, which stubbornly
refused all attempts to study its birational geometry for
about 15 years (there were quite a few well-studied examples of Fano varieties, starting from the quartic, and for the
conic bundles there was Sarkisov's theorem).

The further development of the theory went in two 
directions. Already in the late 80ies V.G.Sarkisov suggested
a general program of factorization of birational maps
between three-dimensional Mori fiber spaces into a composition
of elementary links [S3]. M.Reid did a lot of work to
popularize Sarkisov's ideas among the experts in Mori
theory [R]. In [C1] Corti gave a complete proof of the main
claim of Sarkisov's  program and thus brought the
construction of this theory to an end. Combining the
classical methods with the Sarkisov program, it was proved
in [CPR] that three-dimensional ${\mathbb Q}$-Fano hypersurfaces of
index 1 in weighted projective spaces are birationally 
rigid. Generators of their groups of birational self-maps
were described. Besides, in [C2] Corti suggested to use
the connectedness principle of Shokurov and Koll\' ar
(based, in its turn, on the Kawamata-Viehweg vanishing
theorem) in the investigation of maximal singularities. This
idea turned out to be very fruitful and has already been
used several times, both in dimension three [CM] and
in arbitrary dimension [P12]. This technique was discussed
and further developed in [Ch,ChPk], see also [I2,I3].

The classical technique was developing parallel to the
ideas coming from the Mori theory and the log minimal model
program. In [P7] the construction of hypertangent divisors
was introduced. This construction proved extremely fruitful
and made it possible to prove birational rigidity of generic
Fano varieties and Fano fibrations for several important
families [P8-P11]. This construction makes the basis of the
present paper, either.

This paper was started by the author during his work at the
University of Bayreuth as a Humboldt Research Fellow, and
completed at Max-Planck-Institut f\" ur Mathematik in Bonn.
The author is grateful to Alexander von Humboldt Foundation
and to Max-Planck-Institut f\" ur Mathematik in Bonn for
the financial support of his research. I would like to thank
Mathematisches Institut von Universit\" at Bayreuth and,
in the first place, Prof. Th. Peternell for the warm 
friendly atmosphere which surrounded me in Bayreuth and
for the constant interest to my work. I am very grateful
to Max-Planck-Institut f\" ur Mathematik in Bonn for the
stimulating creative atmosphere, general hospitality and
for exceptionally good conditions of work.


\section{The method of maximal singularities}

\subsection{A criterion of birational superrigidity}

We will prove the main theorem checking the following convenient
sufficient condition of birational superrigidity.
\vspace{0.3cm}

{\bf Proposition 1.1.} {\it Let $X$ be a smooth Fano variety with
$\mathop{\rm Pic} X={\mathbb Z} K_X$. Assume that for any irreducible
subvariety $Y\subset X$ of codimension two
the following two properties are satisfied:} \par
\vspace{0.1cm}

(i) $\mathop{\rm mult}\nolimits_Y\Sigma\leq n$
{\it for any linear system $\Sigma\subset |-nK_X|$
without fixed components};\par
\vspace{0.1cm}

(ii) {\it the inequality}
\begin{equation}
\label{c1}
\mathop{\rm mult}\nolimits_x Y\leq
\frac{\displaystyle 4}{\displaystyle
\mathop{\rm deg} X}
\mathop{\rm deg} Y
\end{equation} 
{\it  holds for any point} $x\in Y$, {\it where 
$$
\mathop{\rm deg} X=(-K_X)^{\mathop{\rm dim} X},\quad
\mathop{\rm deg} Y=(Y\cdot (-K_X)^{\mathop{\rm dim} Y})
$$
and $\mathop{\rm mult}\nolimits_Y\Sigma$ means multiplicity
of a general divisor $D\in\Sigma$ along $Y$. 
\vspace{0.1cm}

Then the variety $X$ is birationally superrigid.} 
\vspace{0.3cm}

It is in this way that birational superrigidity was proved
for the majority of known classes of Fano varieties [P7-P11].

{\bf Proof of Proposition 1.1.} For convenience of the reader
we give it here with all the significant details. For the
other details and comments see [P5,P8].  Assume that
$X$ is not superrigid. Then, by the definition of
superrigidity, we get a birational map
$\chi\colon X-\,-\,\to X'$ and a pair of linear systems
$\Sigma$, $\Sigma'$, transformed one into another by
$\chi$, such that inequality (\ref{i1}) is not true. The next
step in the arguments is given by
\vspace{0.3cm}

{\bf Proposition 1.2.} {\it There exists a geometric
discrete valuation $\nu$ on $X$ such that the {\it Noether-Fano
inequality}
\begin{equation}
\label{c2}
\nu(\Sigma)>n\cdot\mathop{\rm discrepancy}(\nu)
\end{equation}
holds, where $n\in{\mathbb Z}_+$ is defined by the
inclusion $\Sigma\subset |-nK_X|$ and $\nu(\Sigma)=\nu(D)$
for a general divisor $D\in\Sigma$. The discrete
valuations $\nu$, satisfying (\ref{c2}), are called
{\it maximal singularities} of the linear system $\Sigma$.
}
\vspace{0.3cm}

Recall [C2,CPR,P5,P8] that a discerete valuation is said to be 
{\it geometric}, if it is realizable by a prime Weil 
divisor on some model of the field of rational function
${\mathbb C}(X)$.

Let $B\subset X$ be the centre of $\nu$ on $X$. If
$\mathop{\rm codim}\nolimits_XB=2$, then it is easy to see
that $\mathop{\rm mult}\nolimits_B\Sigma>n$, which
contradicts (i). Therefore
$\mathop{\rm codim}\nolimits_XB\geq 3$. Here we come to
the crucial point of the proof.

Let $D_1,D_2\in\Sigma$ be general divisors. They have no
common components and therefore their intersection is of
codimension 2. Denote by $(D_1\circ D_2)$ the effective
cycle of their scheme-theoretic intersection. Obviously,
\begin{equation}
\label{c3}
\mathop{\rm deg}(D_1\circ D_2)=n^2\mathop{\rm deg}X.
\end{equation}

Now the crucial fact is given by
\vspace{0.3cm}

{\bf Proposition 1.3.} {\it The following inequality holds}
\begin{equation}
\label{c4}
\mathop{\rm mult}\nolimits_B (D_1\circ D_2)> 4n^2.
\end{equation}
\vspace{0.3cm}

{\bf Proof} of the proposition is given below. We give here
both known arguments: the classical one, based on the
technique of counting multiplicities, and the recent
argument of Corti, based on his idea of using the
connectedness principle of Shokurov and Koll\' ar.

Now, comparing (\ref{c3}) and (\ref{c4}), we find an irreducible
subvariety $Y\subset X$ of codimension 2 (a component of the
effective cycle $(D_1\circ D_2)$) such that
$$
\mathop{\rm mult}\nolimits_B Y>\frac{4}{\mathop{\rm deg} X}
\mathop{\rm deg} Y,
$$
which contradicts the assumption (ii) of Proposition 1.1.

Therefore, out initial assumption is false and $X$ is superrigid.
Q.E.D. for Proposition 1.1.
 

\subsection{The first proof: counting multiplicities}

First of all, we take a resolution of the discrete valuation
$\nu$. Consider the sequence of blow ups
$$
\begin{array}{cccc}
\displaystyle
\varphi_{i,i-1}: & X_i & \to & X_{i-1} \\
\displaystyle
            & \bigcup & & \bigcup \\
\displaystyle
            & E_i & \to & B_{i-1}
\end{array}
$$
$i\geq 1$, where
$X_0=X, \varphi_{i,i-1}$
blows up the cycle
$B_{i-1}=Z(X_{i-1},\nu)$
of codimension $\geq 2$,
$E_i=\varphi^{-1}_{i,i-1}(B_{i-1})\subset X_i$.
Set also for
$i>j$
$$
\varphi_{i,j}=\varphi_{j+1,j}\circ\dots\circ
\varphi_{i,i-1}:X_i\to X_j,
$$
$$
\varphi_{i,i}=\mathop{\rm id}\nolimits_{X_i}.
$$
For any cycle
$(\dots)$
we denote its proper inverse image on
$X_i$
by adding the upper index $i$:
$(\dots)^i$.
\vspace{0.1cm}

{\bf Remark.}
(i) Note that
$\varphi_{i,j}(B_i)=B_j$
for $i\geq j$.

\noindent
(ii) Note also that although all the $X$'s
are possibly singular,
$B_i\not\subset\mathop{\rm Sing} X_i$
for all $i$.
\vspace{0.1cm}
 
For some $K\in{\mathbb Z}_+$ the divisor
$E_K\subset X_K$ realizes the discrete valuation
$\nu$. 

Let us define the structure of an oriented graph
on the set of exceptional divisors or, equivalently,
on the set of indices $\{1,\dots,K\}$. We draw an
arrow 
$$
i\to j
$$
if $i>j$ and $B_{i-1}\subset E^{i-1}_j$. Set
$p_{ij}$ to be the number of paths from $i$ to
$j$ if $i\neq j$, and $p_{ii}=1$. Finally, set
$p_i=p_{Ki}$
for all
$i=1,\dots,K$.

Now let $\Sigma^j$ be the strict transform of
the linear system $\Sigma$ on $V_j$. Set
$$
\nu_j=\mathop{\rm mult}\nolimits_{B_{j-1}}\Sigma^{j-1}.
$$
Obviously, the multiplicity of the linear system
$\Sigma$ with respect to the valuation $E_j$ is equal to
$$
\nu_{E_j}(\Sigma)=
\sum^j_{i=1}p_{ji}\nu_i.
$$
Setting $\delta_i=\mathop{\rm codim}B_{i-1}-1$,
we get the well-known expression for the discrepancy
$$
a(X,\nu_{E_j})=
\sum^j_{i=1}p_{ji}\delta_i.
$$
The Noether-Fano inequality takes the form
$$
\sum^K_{i=1}p_i\nu_i>
\sum^K_{i=1}p_i\delta_i n.
$$

Now let us consider the following general situation.
Let
$B\subset X, B\not\subset\mathop{\rm Sing} X$
be an irreducible cycle of codimension $\geq 2$,
$\sigma_B\colon X(B)\to X$
be its blowing up,
$E(B)=\sigma^{-1}_B(B)$
the exceptional divisor. Let
$$
Z=\sum m_iZ_i, \quad
Z_i\subset E(B)
$$
be a $k$-cycle,
$k\geq\mathop{\rm dim} B$.
We define the
{\it degree}
of $Z$ as
$$
\mathop{\rm deg} Z=
\sum_im_i\mathop{\rm deg}\left(
Z_i\bigcap\sigma^{-1}_B(b)
\right),
$$
where
$b\in B$
is a generic point,
$\sigma^{-1}_B(b)\cong{\mathbb P}^{\mathop{\rm codim} B-1}$
and the right-hand side degree is the ordinary degree in
the projective space.

Note that
$\mathop{\rm deg} Z_i=0$
if and only if
$\sigma_B(Z_i)$
is a proper closed subset of
$B$.

Our computations will be based upon the following statement.

Let
$D$
and
$Q$
be two different prime Weil divisors on
$X$, $D^B$ and $Q^B$
be their proper inverse images on
$X(B)$.
\vspace{0.1cm}

{\bf Lemma 1.1.}
{\it
{\rm (i)}
Assume that
$\mathop{\rm codim} B\geq 3$.
Then
$$
D^B\circ Q^B=
(D\circ Q)^B+Z,
$$
where
$\circ$
stands for the cycle of the scheme-theoretic intersection,
$\mathop{\rm Supp} Z\subset E(B)$,
and
$$
\mathop{\rm mult}\nolimits_B(D\circ Q)=
(\mathop{\rm mult}\nolimits_BD)
(\mathop{\rm mult}\nolimits_BQ)+\mathop{\rm deg} Z.
$$
{\rm (ii)} Assume that
$\mathop{\rm codim} B=2$. Then
$$
D^B\circ Q^B=Z+Z_1,
$$
where
$\mathop{\rm Supp} Z\subset E(B),
\mathop{\rm Supp}\sigma_B(Z_1)$
does not contain $B$,
and
$$
D\circ Q=
\left[
(\mathop{\rm mult}\nolimits_BD)(\mathop{\rm mult}\nolimits_BQ)+
\mathop{\rm deg} Z
\right]
B+
(\sigma_B)_*Z_1.
$$
}\vspace{0.1cm}

{\bf Proof} follows easily from the standard
intersection theory [Ful].

Now let us come back to our discrete valuation $\nu$.

We divide the resolution
$\varphi_{i,i-1}:X_i\to X_{i-1}$
into
the {\it lower part},
$i=1,\dots,L\leq K$,
when
$\mathop{\rm codim} B_{i-1}\geq 3$,
and the {\it upper part},
$i=L+1,\dots,K$,
when
$\mathop{\rm codim} B_{i-1}=2$.
It may occur that
$L=K$
and the upper part is empty.

Let
$D_1,D_2\in\Sigma$
be two different general divisors. We define a sequence of
codimension 2 cycles on $X_i$'s setting
$$
\begin{array}{l}
\displaystyle
D_1\circ D_2=Z_0,\\
\displaystyle
D^1_1\circ D^2_2=Z^1_0+Z_1,\\
\displaystyle
\vdots\\
\displaystyle
D^i_1\circ D^i_2=
(D^{i-1}_1\circ D^{i-1}_2)^i+Z_i,\\
\displaystyle
\vdots
\end{array}
$$
where
$Z_i\subset E_i$.
Thus for any
$i\leq L$
we get
$$
D^i_1\circ D^i_2=
Z^i_0+Z^i_1+\dots+Z^i_{i-1}+Z_i.
$$
For any
$j>i,j\leq L$
set
$$
m_{i,j}=\mathop{\rm mult}\nolimits_{B_{j-1}}(Z^{j-1}_i)
$$
(the multiplicity of an irreducible cycle along a smaller
cycle is understood in the usual sense; for an arbitrary cycle
we extend the multiplicity by linearity).

Set
$d_i=\mathop{\rm deg} Z_i$.

We get the following system of equalities:
$$
\begin{array}{l}
\displaystyle
\nu^2_1+d_1=m_{0,1},\\
\displaystyle
\nu^2_2+d_2=m_{0,2}+m_{1,2},\\
\displaystyle
\vdots\\
\displaystyle
\nu^2_i+d_i=m_{0,i}+\dots+m_{i-1,i},\\
\displaystyle
\vdots\\
\displaystyle
\nu^2_L+d_L=m_{0,L}+\dots+m_{L-1,L}.
\end{array}
$$
Now
$$
d_L\geq
\sum^K_{i=L+1}
\nu^2_i\mathop{\rm deg}(\varphi_{i-1,L})_*B_{i-1}\geq
\sum^K_{i=L+1}\nu^2_i.
$$

{\bf Definition.}
A function
$a:\{1,\dots,L\}\to{\mathbb R}_+$
is said to be
{\it compatible with the graph structure}, if
$$
a(i)\geq\sum_{j\to i}a(j)
$$
for any
$i=1,\dots,L$.

We will actually use only one compatible function,
namely $a(i)=p_{Ki}=p_i$.
\vspace{0.3cm}

{\bf Proposition 6.}
{\it
Let
$a(\cdot)$
be any compatible function. Then
$$
\sum^L_{i=1}a(i)m_{0,i}\geq
\sum^L_{i=1}a(i)\nu^2_i+
a(L)\sum^K_{i=L+1}\nu^2_i.
$$
}
\vspace{0.3cm}

{\bf Proof.}
Multiply the $i$-th equality by
$a(i)$
and put them all together: in the right-hand part for any
$i\geq 1$
we get the expression
$$
\sum_{j\geq i+1}a(j)m_{i,j}
$$
In the left-hand part for any
$i\geq 1$
we get
$$
a(i)d_i.
$$
\vspace{0.1cm}

{\bf Lemma 1.2.}
{\it
If
$m_{i,j}>0$,
then
$i\to j$.
}
\vspace{0.1cm}

{\bf Proof.} If
$m_{i,j}>0$,
then some component of
$Z^{j-1}_i$
contains
$B_{j-1}$.
But
$Z^{j-1}_i\subset E^{j-1}_i$.
Q.E.D.

Besides, we can compare the
multiplicities
$m_{ij}$
with the corresponding degrees.
\vspace{0.1cm}

{\bf Lemma 1.3.}
{\it
For any
$i\geq 1, j\leq L$
we have
$$
m_{i,j}\leq d_i.
$$
}
\vspace{0.1cm}

{\bf Proof.} The cycles
$B_a$
are non-singular at their generic points. But since
$$
\varphi_{a,b}:B_a\to B_b
$$
is surjective, we can count multiplicities at generic points.
Now the multiplicities are non-increasing with respect to
blowing up of a non-singular cycle, so we are reduced to the
obvious case of a hypersurface in a projective space.
Q.E.D.

Consequently, we have the following estimate:
$$
\sum_{j\geq i+1}a(j)m_{i,j}=
\sum_{
\begin{array}{c}
\scriptstyle
j\geq i+1\\
\scriptstyle
m_{i,j}\neq 0
\end{array}
}a(j)m_{i,j}\leq
d_i\sum_{j\to i}a(j)\leq
a(i)d_i.
$$

So we can throw away all the
$m_{i,*},i\geq 1$,
from the right-hand part, and all the
$d_i, i\geq 1$,
from the left-hand part, replacing
$=$ by $\leq$.
Q.E.D. \vspace{0.3cm}

{\bf Corollary 1.1.}
{\it
Set
$m=m_{0,1}=\mathop{\rm mult}\nolimits_B (D_1\circ D_2), 
D_i\in \Sigma$.
Then
$$
m\left(
\sum^L_{i=1}a(i)
\right)\geq
\sum^L_{i=1}a(i)\nu^2_i+
a(L)\sum^K_{i=L+1}\nu^2_i.
$$
}
\vspace{0.1cm}

{\bf Corollary 2.}
{\it
The following inequality holds:
$$
m\left(
\sum^L_{i=1}p_i
\right)
\geq
\sum^K_{i=1}p_i\nu^2_i.
$$
}
\vspace{0.3cm}
{\bf Proof:}
for
$i\geq L+1$ obviously
$p_i\leq p_L$.
Q.E.D.

Taking into account the Noether-Fano
inequality for
$\nu$,
we see that the right-hand part
here is strictly greater than the
value of the quadratic form
$\sum\limits^K_{i=1}p_i\nu^2_i$
at  the point
$$
\nu_1=\dots=\nu_K=
\frac{
\displaystyle\sum^K_{i=1}p_i\delta_in}
{\displaystyle\sum^K_{i=1}p_i}.
$$
Now set
$$
\Sigma_l=
\sum_{\delta_j\geq 2}p_j,
$$
$$
\Sigma_u=
\sum_{\delta_j=1}p_j.
$$
Then
$$
\mathop{\rm mult}\nolimits_xZ>
\frac{\displaystyle
(2\Sigma_l+\Sigma_u)^2}{\displaystyle
\Sigma_l(\Sigma_l+\Sigma_u)}
n^2.
$$
Now easy computations show that the right-hand
side is not smaller than
$4n^2$.
Q.E.D. for Proposition 2.


\subsection{The second proof: the connectedness 
principle of Shokurov and Koll\' ar}

Here we give an alternative proof of Proposition 2,
suggested by Corti [C2]. Our version is slightly
different from that of [C2], we reduce the proof
to one simple fact about oriented graphs which
was proved originally in [P3].

It is sufficient to consider the case when $B=o$ is a
point. Let $S\ni o$ be a general germ of a smooth
surface on the variety $X$. Obiovuly, $\Lambda=\Sigma|_S$
is a germ of a linear system of curves on $S$ with
the point $o$ as an isolated base point. Since the
pair
$$
(X, \frac{1}{n}\Sigma)
$$
is not canonical at the point $o$ (it is a reformulation
of the Noether-Fano inequality (\ref{c2})), according to
the inversion of adjunction (which is a direct 
consequence of the connectedness principle of
Shokurov and Koll\' ar [K,Sh]), the pair
$$
(S,\frac{1}{n}\Lambda)
$$
is not {\it log}-canonical at the point $o$. In other
words, for a certain birational morphism
$$
\varphi\colon \widetilde S\to S
$$
of smooth surfaces there exists a prime divisor
$E\subset\widetilde S$, satisfying the log-version
of the Noether-Fano inequality
\begin{equation}
\label{c5}
\nu_E(\Lambda)>n(a(E)+1),
\end{equation}
where $a(\cdot)$ is the discrepancy, $\nu_E(\cdot)$
is the multiplicity of a general divisor of the system
at $E$. Let $D_1,D_2\in \Lambda$ be generic curves,
$$
Z=(D_1\circ D_2)
$$
be a zero-dimensional subscheme on $S$. One may assume
that it is supported at the point $o$.
\vspace{0.3cm}

{\bf Proposition 1.4.} {\it The following inequality holds}
$$
\mathop{\rm mult}\nolimits_oZ(=\mathop{\rm deg}\nolimits Z)>4n^2.
$$
\vspace{0.3cm}

Since our considerations are local, 
$\mathop{\rm mult}\nolimits_oZ=\mathop{\rm deg}\nolimits Z$ is just the degree of the zero-dimensional scheme $Z$. Since $S\ni o$ is a general germ of a surface,
we obtain immediately the claim of Proposition 1.3.

{\bf Proof.} We give an elementary argument based on
explicit computations. The original argument of Corti
see in [C2]. Let
$$
\begin{array}{rccc}
   &  S_i &  \supset & E_i  \\
\varphi_{i,i-1} & \downarrow &   &   \downarrow \\
   &  S_{i-1}  & \ni  &  x_{i-1}
\end{array}
$$
be the resolution of the discrete valuation
$\nu_E$, $i=1,\dots,N$, $x_0=o$, $x_1,\dots,x_{N-1}$
points on $S_1,\dots, S_{N-1}$, respectively, where
$$
x_i\in E_i \quad\mbox{and}\quad \nu_{E_N}=\nu_E.
$$
Set $\Gamma$ to be the graph of this resolution:
$$
\{1,\dots,N\}
$$
is the set of vertices, and the vertices $i$ and $j$,
$i>j$, are connected by an oriented edge (notation:
$$
i\to j
$$
always implies that $i>j$), if and only if
$$
x_{i-1}\in E^{i-1}_j,
$$
where $E^{i-1}_j$ is the strict transform of the
exceptional line $E_j$ on $S_{i-1}$. Set also
$$
p_j=(\mbox{the number of paths from}\,\, N\,\,
\mbox{to}\,\, j)
$$
for $j\leq N-1$, $p_N=1$. Set
$$
\nu_i=\mathop{\rm mult}\nolimits_{x_{i-1}}\Lambda^{i-1},
$$
where  $\Lambda^{i-1}$ is the strict transform
of the linear system $\Lambda$ on $S_{i-1}$.
It is easy to see that in terms of the resolution
the log-inequality (\ref{c5}) takes the form
$$
\sum^{N}_{i=1}p_i\nu_i>n\left(
\sum^{N}_{i=1}p_i+1\right).
$$
Besides, the following estimate is true
$$
\mathop{\rm mult}\nolimits_oZ\geq \sum^{N}_{i=1}\nu_i^2.
$$
\vspace{0.1cm}

{\bf Lemma 1.4.} {\it The following inequality holds}
$$
\mathop{\rm mult}\nolimits_oZ>
\frac{\displaystyle
\left(
\sum^{N}_{i=1}p_i+1
\right)^2
}{\displaystyle
\sum^{N}_{i=1}p_i^2
}n^2.
$$
\vspace{0.1cm}

{\bf Proof} is obtained by elementary computations:
the minimum of the quadratic form
$$
\sum^{N}_{i=1}s^2_i
$$
under the restriction
\begin{equation}
\label{c6}
\sum^{N}_{i=1}p_is_i=c
\end{equation}
is attained at $s_i=p_ia$, where the common
constant $a$ is found from (\ref{c6}). Q.E.D. for
the lemma.

In view of this lemma Proposition 1.4 is an implication
of the following fact.
\vspace{0.1cm}

{\bf Lemma 1.5.} {\it The following inequality holds}
\begin{equation}
\label{c7}
\left(\sum^{N}_{i=1}p_i+1\right)^2\geq
4\sum^{N}_{i=1}p_i^2.
\end{equation}
\vspace{0.1cm}

{\bf Proof.} Note first of all that in (\ref{c7})
the equality can be attained, for instance when
$N=1$. Assume that $N\geq 2$. Set
$$
\{2,\dots,k\leq N\}=\{i\,|\,i\to 1\}.
$$
By the definition of the integers $p_i$ we get
the equality
$$
p_1=\sum_{i\to1}p_i=\sum^{k}_{i=2}p_i.
$$
Consequently, (\ref{c7}) can be rewritten as
$$
\left(
2p_1+\sum^{N}_{i=k+1}p_i+1\right)^2\geq
4\sum^{N}_{i=1}p_i^2
$$
or, after an easy transformation,
$$
4\left(
\sum^{k}_{i=2}p_i\right)
\left(\sum^{N}_{i=k+1}p_i+1\right)+
\left(\sum^{N}_{i=k+1}p_i+1\right)^2\geq
4\sum^{k}_{i=2}p_i^2+4\sum^{N}_{i=k+1}p^2_i.
$$
It is easy to see that if $k=N$, then the subgraph
of the graph $\Gamma$ with the vertices 
$\{2,\dots,N\}$ is of the form
$$
2\leftarrow 3\leftarrow \dots\leftarrow N
$$
(since on any surface $S_i$ the curve
$$
\mathop{\bigcup}\limits_{j\leq i} E^i_j
$$
is by smoothness a normal crossing divisor).
Hence $p_2=\dots=p_N=1$ and the inequality (\ref{c7})
holds in an obvious way. So let us assume that
$N\geq k+1$. Arguing by induction on the number
of vertices of the graph $\Gamma$ we may assume
that the inequality
$$
\left(
\sum^{N}_{i=k+1}p_i+1\right)^2\geq
4\sum^{N}_{i=k+1}p_i^2
$$
is true. Therefore, it is enough to show that
the following estimate is true:
\begin{equation}
\label{c8}
\left(\sum^{k}_{i=2}p_i\right)
\left(\sum^{N}_{i=k+1}p_i+1\right)\geq
\sum^{k}_{i=2}p_i^2.
\end{equation}

If $k=2$, then by construction we get
$$
p_2\leq\sum^{N}_{i=3}p_i,
$$
which immediately implies (\ref{c8}). If $k\geq 3$,
then the subgraph of the graph $\Gamma$ with the
vertices $\{2,\dots,\}$ is a chain:
$$
2\leftarrow 3\leftarrow \dots\leftarrow k.
$$
Since $k\to (k-1)$ and $k\to 1$, the vertices $k$
and $i$, $i\in\{2,\dots,k-2\}$, are not joined by
an arrow (oriented edge). Consequently,
$$
j\not\to i
$$
for $j\geq k+1$,  $i\in\{2,\dots,k-2\}$. Thus each
path from the vertex $N$ to the vertex
$i\in\{2,\dots,k-2\}$  must go through the vertex
$k-1$. Therefore
\begin{equation}
\label{c9}
p_2=\dots=p_{k-1}=p_k+\sum_{
\begin{array}{c}\scriptstyle
i\to k-1 \\  \scriptstyle
i\geq k+1
\end{array}
}p_i.
\end{equation}

{\bf Lemma 1.6.} {\it  For any $i\in\{1,\dots,N\}$ the
following inequality holds:
\begin{equation}
\label{c10}
p_i\leq\sum_{j\geq i+2}p_j+1
\end{equation}
(if the set $\{j\geq i+2\}$ is empty, the
sum is assumed to be equal to zero).}
\vspace{0.1cm}

{\bf Proof} is obtained by decreasing induction on $i$.
If $i=N$ or $i=N-1$, then the inequality (\ref{c10}) is
true. Now we get
$$
p_i-\sum_{j\geq i+2}p_j=
\sum_{j\to i}p_j
- \sum_{j\geq i+2}p_j=
$$
$$
=p_{i+1}-\sum_{
\begin{array}{c} \scriptstyle
j\geq i+2 \\  \scriptstyle
j\not\to i
\end{array}
}p_j.
$$
Write down the set $\{j\,|\,j\to i\}$ as
$\{i+1,\dots,i+l\}$. If $l=1$,
then applying the induction hypothesis, we obtain
(\ref{c10}). If $l\geq 2$, then similarly to (\ref{c9}) we
get
$$
p_{i+1}=\dots=p_{i+l-1}=p_{i+l}+\sum_{
\begin{array}{c}  \scriptstyle
j\to i+l-1 \\  \scriptstyle
j\geq i+l+1
\end{array}
}p_j.
$$
Therefore
$$
p_{i+1}-\sum_{
\begin{array}{c}  \scriptstyle
j\geq i+2 \\  \scriptstyle
j\not\to i
\end{array}
}p_j=p_{i+l-1}-
\sum_{j=i+l+1}^Np_j.
$$
Applying the induction hypothesis to $i+l-1$, we
complete the proof.

Let us come back to the proof of Lemma 1.5. We get
$$
p_2=\dots=p_{k-1}\leq \sum_{i=k+1}^N p_i+1.
$$
But $p_k\leq p_{k-1}$, so that
$$
\left(
\sum^k_{i=2}p_i\right)\left(
\sum^N_{i=k+1}p_i+1\right)\geq
p_{k-1}\sum^k_{i=2}p_i\geq \sum^k_{i=2}p^2_i,
$$
which is what we need.

The proof of Proposition 1.4 is complete.

{\bf Remark.} The inequality (\ref{c10}) first
appeared in [P3] as an auxiliary estimate. Here we
have reproduced the inductive proof of this fact
given in [P3] for convenience of the reader.


\subsection{The Lefschetz theorem}

Let us prove that the condition (i) of Proposition 2
is true for any smooth iterated double cover $V$.
Indeed, let $\Sigma\subset|-nK_V|$ be a linear system of divisor
without fixed components. Let $Y\subset \mathop{\rm Bs}\Sigma$ be
an irreducible subvariety of codimension 2. Take two general divisors,
$D_1,D_2\in\Sigma$ and consider the algebraic cycle of their
scheme-theoretic intersection:
$$
Z=(D_1\circ D_2)=aY+\sum a_iY_i,
$$
where $a_i\geq 1$, $Y_i\neq Y$ are some irreducible subvarieties of
codimension 2. Obviously,
$a\geq (\mathop{\rm mult}\nolimits_Y\Sigma)^2$. Compute the degrees:
$$n^2\mathop{\rm deg} V=\mathop{\rm deg} Z=
a\mathop{\rm deg} Y+\sum a_i\mathop{\rm deg} Y_i,$$
where $\mathop{\rm deg} V=2^m d_1\dots d_k$. By the Lefschetz theorem
$Y$ is numerically equivalent to $mK^2_V$, $m\geq 1$, so 
that
$\mathop{\rm deg} Y=m\mathop{\rm deg} V$, whence $$
(\mathop{\rm mult}\nolimits_Y\Sigma)^2\leq a \leq n^2/m\leq n^2,$$
which is what we need. Q.E.D. for the condition (i).


\subsection{Proof of Proposition 0.1}

Part (i) of Proposition 0.1 is almost obvious.
If
$$
\begin{array}{ccccl}
\chi\colon & X & -\,-\,\to & X' &  \\
           &   &           & \downarrow & \pi' \\
           &   &           &  S' &
\end{array}
$$
is a birational map onto $X'$, where $\dim S'\geq 1$
and fibers of $\pi'$ are uniruled, take $\Sigma'$ to be
a pull back of a moving linear system on $S'$. Then
$c(\Sigma',X')=0$ and therefore by superrigidity
$c(\Sigma,X)=0$. But $\Sigma\subset |-nK_X|$ is a
moving linear system, hence $c(\Sigma,X)=n\geq 1$.
Contradiction.

Let us prove part (ii) of Proposition 0.1. 
Let $\chi\colon X-\,-\,\to X'$ be a birational map,
$\varphi\colon Y\to X$ be its Hironaka resolution, so that
$\psi=\chi\circ\varphi\colon Y\to X'$ is a birational morphism.
The variety $Y$ is non-singular and
$$
\mathop{\rm Pic}Y={\mathbb Z}\varphi^*K_X\oplus
\mathop{\bigoplus}\limits_{i\in I}{\mathbb Z} E_i,
$$
where $\{E_i|i\in I\}$ is the set of all the
$\varphi$-exceptional divisors. By assumption,
$$
\mathop{\rm Pic}Y\otimes{\mathbb Q}={\mathbb Q}\psi^*K_{X'}\oplus
\mathop{\bigoplus}\limits_{j\in J}{\mathbb Q} E'_j,
$$
where $\{E'_j|j\in J\}$ is the set of all the
$\psi$-exceptional divisors. For simplicity of notations
set $K=\varphi^* K_X$, $K'=\psi^* K_{X'}$. We get
\begin{equation}
\label{c11}
K_Y=K+\sum_{i\in I}a_iE_i=K'+\sum_{j\in J}a'_jE'_j,
\end{equation}
where $a_i\in{\mathbb Z}$, $a_i\geq 1$, and $a'_j\in{\mathbb Q}$,
$a'_j>0$. Let $\Sigma'=|-mK_{X'}|$, $m\gg 0$, be a very
ample system. Obviously, $c(\Sigma',X')=m$. Take
$\Sigma=\chi^{-1}_*\Sigma'\subset|-nK_X|$; obviously,
$c(\Sigma,X)=n$. Twisting by a suitable birational 
self-map, we may assume that the inequality (\ref{i1})
is already satisfied for $\chi$. Hence $n\leq m$. The proper inverse
image of the linear system $\Sigma$ on $Y$ coincides with
the inverse image of the linear system $\Sigma'$ with
respect to $\psi$. Therefore, there exist positive integers
$b_i$, $i\in I$, such that
$$
-mK'=-nK-\sum_{i\in I}b_iE_i.
$$
Dividing by $m$ and substituting into (\ref{c11}), we get
$$
\left(1-\frac{n}{m}\right)K=\sum_{i\in I}\left(
\frac{b_i}{m}-a_i\right)E_i+\sum_{j\in J}a'_jE'_j.
$$
Since the divisors $E_i$ are $\varphi$-exceptional and
$a'_j>0$, we get the equality $n=m$: otherwise we get a
contradiction with the ampleness of $(-K_X)$. Furthermore,
all the divisors $E'_j$ turn out to be $\varphi$-exceptional
and, moreover, $\{E_i|i\in I\}=\{E'_j|j\in J\}$, otherwise
$\mathop{\rm rk} \mathop{\rm Pic} X'\geq 2$. Thus $\chi$
is an isomorphism in codimension one: set
$$
U=X\setminus \mathop{\bigcup}\limits_{i\in I}\varphi(E_i),
\quad
U'=X'\setminus \mathop{\bigcup}\limits_{j\in J}\psi(E'_j),
$$
then $\chi\colon U\to U'$ is an isomorphism. Therefore
$\Sigma=|-nK_X|$ and $\chi$ induces an isomorphism of the
linear systems $\Sigma$ and $\Sigma'$. Consequently,
$\chi\colon X\to X'$ is an isomorphism. (Strictly speaking,
we have proved that for an arbitrary birational map
$\chi\colon X-\,-\,\to X'$ there exists 
$\chi^*\in\mathop{\rm Bir} X$ such that 
$\chi\circ\chi^*$ is an isomorphism.) The rest is
obvious. Proof of Proposition 0.1 is complete.

Proof of Proposition 0.2 is similar to the elementary
arguments for the part (i) of Proposition 0.1 above.


\section{Iterated double covers}

\subsection{Coordinate presentations}

Let ${\mathbb P}^{\sharp}$ be the weighted projective space
$$
{\mathbb P}(\underbrace{1,1,\dots,1}_{M+k+1},l_1,\dots,l_m),
$$
where to the weights $l_i\geq 2$ correspond the new
homogeneous coordinates $y_i$, $i=1,\dots,m$. The variety
${\mathbb P}^{(m)}$ can be realized as a complete intersection
of the type $2l_1\cdot\dots\cdot 2l_m$ in ${\mathbb P}^{\sharp}$:
$$
{\mathbb P}^{(m)}=\mathop{\bigcap}\limits^m_{i=1}
\{y_i^2=g_i\}\subset {\mathbb P}^{\sharp}.
$$
In a similar way, the variety 
$V\subset{\mathbb P}^{(m)}\subset {\mathbb P}^{\sharp}$ is a complete 
intersection of the type
$d_1\cdot\dots\cdot d_k\cdot 2l_1\cdot\dots\cdot 2l_m$.

Let $p\in Q$ be an arbitrary point,
$(z_1,\dots,z_{M+k})$ a system of affine coordinates with
the origin at the point $p$. Set
$$
\begin{array}{l}
f_i=g_{i,1}+\dots+g_{i,d_i}, \\   \\  g_i=w_{i,0}+w_{i,1}+\dots+w_{i,2l_i}
\end{array}
$$
to be the Taylor decompositions of the (non-homogeneous)
polynomials $f_i$, $g_i$ in the coordinates $z_*$. Here
$$
\mathop{\rm deg}\nolimits q_{i,j}=j,\quad \mathop{\rm deg}\nolimits w_{i,j}=j.
$$
If $w_{i,0}\neq 0$, that is, $p\not\in W_i$, then for the
convenience of computations we assume always that
$w_{i,0}=1$.

{\bf Definition 2.1.} The point $p\in Q$ is of class 
$e\in\{0,1,\dots, m\}$, if
$$
e=\sharp\{i\,|\, p\in W_i\}.
$$
We write the set ${\cal L} ={\cal L}(p)=\{i\,|\, p\in W_i\}$
as
$$
\{i_1<\dots<i_e\}.
$$
Let us define a convenient coordinate system at a point
$p\in Q$ of class $e$. For simplicity assume that
$z_i=x_i/x_0$, $i=1,\dots,M+k$. Set
$$
u_i=y_i/x_0^{l_i}
$$
for $i=1,\dots,m$. The set of regular functions $(z_*,u_*)$
is a system of affine coordinates on an open affine subset
$U\subset{\mathbb P}^{\sharp}$, $U\cong{\mathbb C}^{M+k+m}$.
With respect to these coordinates the variety $V$ is given
by the system of equations
\begin{equation}
\label{a1}
\left\{
\begin{array}{ll}
f_i(1,z_1,\dots,z_{M+k})=0, & i=1,\dots,k, \\  \\
u_i^2=g_i(1,z_1,\dots,z_{M+k}), & i=1,\dots,m,
\end{array}
\right.
\end{equation}
and from now on we will identify the homogeneous
polynomials $f_i$, $g_i$ with their non-homogeneous
presentations of the type $f(1,z_*)$.

If the point $p\in Q$ is of class $0$, then for all
$i=1,\dots,m$ we have $g_i(p)\neq 0$, that is, 
$w_{i,0}=1$. In this case for any point 
$q\in\sigma^{-1}(p)$ the linear maps
$$
\sigma_*\colon T_q{\mathbb P}^{(m)}\to T_p{\mathbb P},
$$
$$
\sigma_*\colon T_q V\to T_pQ
$$
are isomorphisms, so that the $\sigma$-preimage of
any system of local coordinates on ${\mathbb P}$ and $Q$
makes a system of local coordinates on ${\mathbb P}^{(m)}$
(for instance, $(z_1,\dots,z_{M+k})$)
and $V$, respectively.

If the point $p\in Q$ is of class $e\geq 1$, assume
that ${\cal L}(p)=\{1,\dots,e\}$. In this case a
natural system of local coordinates on ${\mathbb P}^{(m)}$
is given by the set of functions
$$
(z_{j_1},\dots, z_{j_{M+k-e}},u_1,\dots,u_e),
$$
where $z_{j_1},\dots, z_{j_{M+k-e}}$ make a system
of local coordinates on the complete intersection
$W_1\cap\dots\cap W_e\subset{\mathbb P}$.
\vspace{0.1cm}

{\bf Lemma 2.1.} {\it For
$i\in \{1,\dots,e\}$ the inverse image of the hyperplane $\{w_{i,1}=0\}$ is tangent to $V$ at the point 
$q\in\sigma^{-1}(p)$. In particular, the following
isomorphism holds
$$
T_qV\cong T_p(Q\cap W_1\cap \dots \cap W_e)
\oplus \langle u_1,\dots,u_e\rangle^*.
$$
With respect to this isomorphism the tangent cone to
the intersection $V\cap \{w_{i,1}=0\}$ is given by
the quadratic equation}
$$
u_i^2=w_{i,2}|_{T_p(Q\cap W_1\cap\dots\cap W_e)}.
$$

{\bf Proof.} It is obvious from the system of equations
(\ref{a1}).


\subsection{The regularity condition}

Let $g(z_*)=1+w_1+\dots+w_{2l}$ be a polynomial. Following
[P8], consider the formal series
$$
(1+t)^{1/2}=
1+\sum^{\infty}_{i=1}\gamma_i t^i=
1+\frac12 t-\frac18 t^2+\dots,
$$
and make the following formal series in the variables $z_*$:
$$
\sqrt{g}=
(1+w_1+\dots+w_{2l})^{1/2}=
1+\sum^{\infty}_{i=1}\gamma_i (w_1+\dots+w_{2l})^i=
$$
$$
=1+\sum^{\infty}_{i=1}\Phi_i(w_1,\dots,w_{2l}),
$$
where $\Phi_i(w_1(z_*),\dots,w_{2l}(z_*))$ are
homogeneous polynomials of degree $i$ in $z_*$.
Obviously,
$$
\Phi_i(w_*)=\frac12 w_i+
(\mbox{polynomials in}\,\,\, w_{1},\dots,w_{i-1}).
$$
For instance, $\Phi_1(w_*)=\frac12 w_1$. Furthermore,
for $j\geq 1$ set
$$
[\sqrt{g}]_j=1+\sum^j_{i=1}\Phi_i(w_*(z_*))
$$
and
$$
g^{(j)}=g-[\sqrt{g}]_j^2.
$$

It is easy to see that the first non-zero homogeneous component of $g^{(j)}$ is of degree $j+1$. Denote it
by the symbol $h_{j+1}[g]$. Obviously,
\begin{equation}
\label{a2}
h_{j+1}[g]=w_{j+1}+A_j(w_1,\dots,w_j),
\end{equation}
where we are not interested in the particular structure
of the polynomial $A_j$.

Now let us formulate the regularity condition. Set
$h_{i,j}=h_j[g_i]$.

{\bf Definition 2.2.} (i) A point $p\in Q$ of class $e=0$
is regular with respect to the set $(f_*;g_*)$ if the
set of homogeneous polynomials
\begin{equation}
\label{a3}
\{q_{i,j}\,|\, (i,j)\in{\cal J}_q\}\cup
\{h_{i,j}\,|\, (i,j)\in{\cal J}_h\},
\end{equation}
where
$$
\begin{array}{l}
{\cal J}_q=\{(i,j)\,|\, 1\leq i\leq k, 1\leq j\leq d_i\}, \\
\\ {\cal J}_h=\{(i,j)\,|\,
1\leq i\leq m, l_i+1\leq j\leq 2l_i, (i,j)\neq (m,2l_m)\},
\end{array}
$$
is regular at $p=o=(0,\dots,0)\in{\mathbb C}^{M+k}$,
that is, makes a regular sequence in ${\cal O}_{p,{\mathbb P}}$ or,
in other words, the set of its common zeros is 
one-dimensional.

(ii) A point $p\in Q$ of class $e\geq 1$ is regular with
respect to the set $(f_*;g_*)$, if the set of homogeneous
polynomials
\begin{equation}
\label{a4}
\{q_{i,j}\,|\, (i,j)\in{\cal J}_q\}\cup
\{w_{i,1}\,|\, 1\leq i\leq e\}\cup
\{h_{i,j}\,|\, (i,j)\in{\cal J}_e\},
\end{equation}
 makes a regular sequence in ${\cal O}_{p,{\mathbb P}}$. Here
for simplicity of notations we assume that
$$
{\cal L}(p)=\{1,\dots,e\},
$$
the set ${\cal J}_q$ was defined above and
$$
{\cal J}_e=\{(i,j)\,|\,
e+1\leq i\leq m, l_i+1\leq j\leq 2l_i\}.
$$
\vspace{0.3cm}

{\bf Proposition 2.1.} {\it For a general set 
$(f_*;g_*)\in{\cal H}$ any point $p\in Q(f_*)$ is regular.}
\vspace{0.3cm}

{\bf Proof.} Let us consider, to begin with, the general
problem of estimating the codimension of ``incorrect''
sets of polynomials. Here we follow [P10]. This problem
is of an independent interest.


\subsection{Non-regular sets of polynomials}

Let $z_1,\dots,z_{N+1}$ be a set of variables. The symbol
${\cal P}_a$ stands for the space of homogeneous polynomials
of degree $a$ in the variables $z_*$. Set
$${\cal L}=\prod\limits^{l+1}_{i=1}{\cal P}_{m_i}=
\{(p_1,\dots,p_{l+1})\}$$
to be the set of all $(l+1)$-uples of homogeneous
polynomials in the variables $z_*$,  $0\leq l\leq N-1$.
With each $(l+1)$-uple $(p_*)\in{\cal L}$ we
associate the projectivized set of its zeros
$$Z(p_*)=\{p_1=\dots=p_{l+1}=0\}\subset
{\mathbb P}^N={\mathbb P}({\mathbb C}^{N+1})=X.$$
Here we write $Z(p_*)$ and not $V(p_*)$ in order to
make this notation different from our complete intersection
$V=V(f_*)$, the principal object of study in this paper. Let
$$
Y=\{(p_*)\in{\cal L}|\mathop{\rm codim}\nolimits_XZ(p_*)\leq l\}
$$
be the set of ``irregular'' $(l+1)$-uples. We need an estimate for
the codimension of $Y$. The case $l=N-1$, when the ``correct''
dimension of $Z(p_*)$ is zero, is especially important for
applications to Fano varieties. However, for technical reasons,
it is more convenient to consider the general case of an
arbitrary $l\in\{0,\dots,N-1\}$.
Set $I=\{1,\dots,l+1\}$ and
$$
\mu_j=\mathop{\rm min}\limits_{S\subset I, \sharp S=j}
\left\{\sum_{i\in S}m_i\right\},
$$
$b=1,\dots,l+1$. Assume that
$$m=\mu_1=\min\{m_1,\dots,m_{l+1}\}\geq 2.$$
\vspace{0.3cm}

{\bf Proposition 2.2.} {\it For any $l\in\{1,\dots,N-1\}$
the following estimate holds:}
\begin{equation}
\label{a5}
\mathop{\rm codim}\nolimits_{{\cal L}}Y\geq
\mathop{\rm min}\limits_{j\in\{0,\dots,l\}}
\{(\mu_{j+1}-j)(N-j)+1\}.
\end{equation}
\vspace{0.3cm}

{\bf Remark.} (i) For $l=0$ we have the trivial estimate
$\mathop{\rm codim}\nolimits_{{\cal L}}Y=\mathop{\rm dim}
{\cal P}_m$.

\noindent
(ii) It is not so easy to say whether the estimate (\ref{a5})
is optimal or not. It is obtained below by the method which
is completely different from the technique used in [P7] for a
similar purpose, that is, to prove existence of regular
hypersurfaces $V_M\subset {\mathbb P}$. The technique of
[P7] does not work here: the resulting estimates are too weak
for complete intersections; however, it is possible that
a combination of the method of this paper and that of [P7]
could improve (\ref{a5}).

{\bf Proof of Proposition 2.2.} Set
$${\cal L}_a=\prod\limits^{a}_{i=1}{\cal P}_{m_i}=
\{(p_1,\dots,p_a)\}.$$
For each irregular $(l+1)$-uple $(p_1,\dots,p_{l+1})$ fix the
first (counting from the left to the right) moment when the
codimension of the set of zeros $p_1=\dots=p_a=0$ fails to take the
correct value. Consider the sets
$$Y_a=\{(p_*)\in{\cal L}_a|
\mathop{\rm codim}\nolimits_XZ(p_1,\dots,p_a)=
\mathop{\rm codim}\nolimits_XZ(p_1,\dots,p_{a-1})=a-1\}.$$
Obviously,
$$
Y=\coprod^{l+1}_{a=1}
\left(Y_a\times \prod\limits^{l+1}_{i=a+1}{\cal P}_{m_i}\right).
$$
In particular,
$$\mathop{\rm codim}\nolimits_{{\cal L}}Y=
\min\{\mathop{\rm codim}\nolimits_{{\cal L}_a}Y_a|1\leq a\leq l+1\}.
$$ 
Set $I_a=\{1,\dots,a\}\subset I$ and
$$
\mu_{a,j}=\mathop{\rm min}\limits_{S\subset I_a, \sharp S=j}
\left\{\sum_{i\in S}m_i\right\},
$$
$j=1,\dots,a$. Obviously, $\mu_{a,j}\geq\mu_j$. Therefore, it is
sufficient to prove the estimate
\begin{equation}
\label{a6}
\mathop{\rm codim}\nolimits_{{\cal L}_a}Y_a\geq
\mathop{\rm min}\limits_{j\in\{0,\dots,a-1\}}
\{(\mu_{a,j+1}-j)(N-j)+1\}
\end{equation}
for each $a=2,\dots,l+1$. We omit the trivial case $a=1$,
because in this case
$\mathop{\rm codim}\nolimits_{{\cal L}_1}Y_1=
\mathop{\rm dim}{\cal P}_{m_1}\geq\mathop{\rm dim}{\cal P}_{m}$,
which is certainly higher than the right-hand side of
(\ref{a5}), just set in (\ref{a5}) $j=0$.

The space ${\cal L}_a$, the set $Y_a$ and the inequality
(\ref{a6}) do not depend on $l$. Thus we may simplify our
notations, setting $a=l+1$ and $\mu_{a,j}=\mu_j$. In other
words, we prove the inequality (\ref{a5}) for $Y_{l+1}$
instead of $Y$. 

Denote $Y_{l+1}$ by $Y^*$. We have reduced our original
problem to a simpler task of estimating codimension of $Y^*$
in ${\cal L}$, where $Y^*$ consists of all such $(l+1)$-uples
of polynomials $(p_1,\dots,p_{l+1})$ that the set
$p_1=\dots=p_l=0$ has the correct dimension and there exists
an irreducible component $B\subset Z(p_1,\dots,p_l)$,
on which $p_{l+1}$ vanishes. Let $\langle B\rangle$ be the
linear span of $B$, and set
$\mathop{\rm codim}\nolimits\langle B\rangle=b\leq l$.

Now set $Y^*(b)$ to be the set of all those $(l+1)$-uples
$(p_*)\in Y^*$, for which there exists a component
$B\subset Z(p_1,\dots,p_l)$ such that
$$\mathop{\rm codim}\nolimits_X\langle B\rangle=b, \quad
p_{l+1}|_B\equiv 0.$$
Obviously,

$$
Y^*=\mathop{\bigcup}\limits^l_{b=0}Y^*(b).
$$
Thus it is sufficient to prove that
\begin{equation}
\label{a7}
\mathop{\rm codim}\nolimits_{{\cal L}}Y^*(b)\geq
(\mu_{b+1}-b)(N-b)+1.
\end{equation}
for each $b=0,\dots,l$. Let us prove (\ref{a7}).
\vspace{0.3cm}

{\bf The case ${\bf b=0}$.} Here $\langle B\rangle={\mathbb P}^N$
and
therefore each non-zero monomial in the linear forms in
$z_1,\dots,z_{N+1}$ of degree $d_{l+1}\geq d$ does not vanish
on $B$. The space of monomials $$\left\{
\prod^{m_{l+1}}_{i=1}(a_{i,1}z_1+\dots+a_{i,N+1}z_{N+1})
\right\}\subset{\cal P}_{m_{l+1}}$$
is closed. Its dimension is equal to
$$
m_{k+1}N+1\geq \mu_1N+1.
$$
On the other hand, the set of polynomials
$p_{l+1}\in{\cal P}_{m_{l+1}}$ that vanish on $B$ is closed.
These two closed sets intersect each other at zero only.
Therefore the codimension of $Y(0)$ in ${\cal L}$ is no
smaller than $\mu_1N+1$. This gives (\ref{a7}) for $b=0$.
\vspace{0.3cm}

{\bf The case ${\bf b\geq 1}$.} Here $\langle B\rangle={\mathbb P}^{N-b}$.
Our strategy is to reduce this case to the previous one ($b=0$),
restricting the polynomials $p_i$ to the linear span
$P=\langle B\rangle$. Although our arguments are rather
simple,
they are not straightforward and require some extra work.

{\bf Definition 3.} Let $g_1,\dots,g_e$ be homogeneous polynomials
on the projective space $P$, $e\leq\mathop{\rm dim} P-1$,
$\mathop{\rm deg} g_i\geq 2$ for $i=1,\dots,e$. An irreducible
subvariety $C\subset P$ such that $\langle C\rangle=P$ and
$\mathop{\rm codim}\nolimits_P C=e$ is called an
{\it associated subvariety of the sequence} $(g_*)$, if
there exists a chain of irreducible subvarieties
$R_j\subset P$, $j=0,\dots,e$, satisfying the following
properties:

\begin{itemize}
\item $R_0=P$;
\item for each $j=0,\dots,e-1$ the subvariety $R_{j+1}$
is an irreducible component of the closed algebraic set
$$
\{p_{j+1}=0\}\cap R_j,
$$
where $p_{j+1}|_{R_j}\not\equiv 0$, so that
$\mathop{\rm codim}\nolimits_P R_j=j$ for all $j$;
\item  $R_e=C$.
\end{itemize}

If the sequence $(g_*)$ has an associated subvariety,
this sequence is said to be {\it good}.
\vspace{0.1cm}

{\bf Lemma 2.2.} (i) {\it The property of being good is an
open property.} \par
(ii) {\it A good sequence $(g_*)$ can have at most
$$
\left[\frac{1}{e+1}\prod^e_{j=1}\mathop{\rm deg} g_j
\right]
$$
associated subvarieties.}
\vspace{0.1cm}

{\bf Proof} is easily obtained by induction on $e$. For
$g_1$ we have the condition $g_1\not\equiv 0$, which is clearly
an open one. Furthermore, at least one irreducible
component of the hypersurface $g_1=0$ must be of degree
$\geq 2$, which is also an open condition. There can be
at most $[\mathop{\rm deg}g_1/2]$ such components.

Assume that Lemma is true for each $e=1,\dots,j$, where
$j\leq\mathop{\rm dim} P-2$. Denote by $G_j$ the open
set of good sequences of length $j$. By (ii), for each
$(g_*)\in G_j$ there exist at most
$$
\left[\frac{1}{j+1}\prod^j_{\alpha=1}\mathop{\rm deg} g_{\alpha}
\right]
$$
associated subvarieties. The polynomial $g_{j+1}$ should be
non-zero on at least one of them, say $R_j$, and moreover, the
intersection
$$
\{g_{j+1}=0\}\cap R_j
$$
should contain an irreducible component, the linear span of
which is $P$. Obviously, this determines an open set in
$$
G_j\times H^0(P,{\cal O}_P(\mathop{\rm deg} g_{j+1})).
$$
Each associated subvariety has codimension $j+1$ and does not
lie in a hyperplane; therefore, its degree is not smaller
than $j+2$. Q.E.D. for the lemma.

Now let us return to the polynomials $p_*$ and assume that
$l>b$. We claim that we can find $(l-b)$ polynomials among
them --- after re-numbering we may assume that they are
$p_1,\dots,p_{l-b}$ --- such that the sequence
\begin{equation}
\label{a8}
p_1|_P,\dots,p_{l-b}|_P
\end{equation}
is good and $B$ is one of its associated subvarieties.
\vspace{0.3cm}

{\bf Proof of the claim.} Assume that we have already found
$j$ polynomials --- let them be
$p_1,\dots,p_{j}$ --- such that the sequence
$(p_1|_P,\dots,p_j|_P)$ is good and one of its associated
subvarieties, say $R_j$, contains $B$. If $j<l-b$, then
$R_j\neq B$ and there exists a polynomial $p_{\alpha}$,
$\alpha\in\{j+1,\dots,l\}$, such that
$$
p_{\alpha}|_{R_j}\not\equiv 0.
$$
Otherwise, $R_j\subset Z(p_1,\dots,p_l)$ and we get a
contradiction, since $\mathop{\rm dim}R_j>\mathop{\rm dim}B$.
After re-numbering, we may assume that $\alpha=j+1$. Now
$p_{j+1}|_B\equiv 0$, so that there exists an irreducible
component $R_{j+1}$ of the set $\{p_{j+1}=0\}\cap R_j$,
such that $R_{j+1}\supset B$. Proceeding in this way, we
obtain our claim.

Now fix a projective subspace $P\subset{\mathbb P}^N$
of codimension $b$. Let $Y^*(P)$ be the set of all
$(l+1)$-uples $(p_1,\dots,p_{l+1})\in Y^*$
such that there exists a component
$B\subset Z(p_1,\dots,p_l)$, whose linear span is
$\langle B\rangle = P$ and $p_{l+1}|_B\equiv 0$.

By Lemma 2.2, good sequences form an open set. Thus we
may estimate the codimension of $Y^*(P)$ in ${\cal L}$,
assuming that $(p_1|_P,\dots,p_{l-b}|_P)$ make a good
sequence. Let
$$
B_1,\dots,B_K
$$
be all its associated subvarieties, whose linear span
is $P$. If $(p_1,\dots,p_{l+1})\in Y^*(P)$, then the
polynomials
$$
p_{l-b+1}|_{P},\dots,p_{l+1}|_{P}
$$
must all vanish on one of these subvarieties $B_i$. Now
arguing as in the case $b=0$, we get
$$
N\sum^{l+1}_{j=l-b+1}\mathop{\rm deg} p_j+b+1\geq
\mu_{b+1}(N-b)+b+1
$$
independent conditions on $p_{l-b+1},\dots,p_{l+1}$. Taking
into account that the Grassmanian has dimension
$\mathop{\rm dim} G(N+1-b,N+1)=b(N+1-b)$,
we get finally
$$
\begin{array}{ccc}
\mathop{\rm codim}\nolimits_{\cal L} Y^*(b) & \geq  &
\mu_{b+1}(N-b)+b+1-b(N+1-b)\\ \\
   &   & \parallel \\  \\ &   & (\mu_{b+1}-b)(N-b)+1,
\end{array}
$$
which is what we need.

In our arguments above we assumed that $l>b$. If $l=b$, then
$B\subset{\mathbb P}^N$ is a line, $l=N-1$ and the inequality
(\ref{a7}) can be obtained by an easy dimension count: for
a fixed line $B$ the
condition $p|_B\equiv 0$ for a polynomial $p$ of degree
$e\geq 1$ defines a closed algebraic set of polynomials of
codimension $e+1$ in ${\cal P}_e$. Therefore,
$$
\mathop{\rm codim}\nolimits_{{\cal L}}Y^*(N-1)\geq
\sum^N_{i=1}(m_i+1)-2(N-1)=\mu_{l+1}-l+1,
$$
since $\mu_{l+1}=m_1+\dots+m_N$.

Q.E.D. for Proposition 2.2.
\vspace{0.3cm}

{\bf Corollary 2.1.} {\it In the notations of Proposition 2.2
for $l\leq N-2$ the following estimate holds:
\begin{equation}
\label{a9}
\mathop{\rm codim}\nolimits_{\cal L} Y\geq
mN+1,
\end{equation}
whereas for $l=N-1$ the following estimate holds:}
\begin{equation}
\label{a10}
\mathop{\rm codim}\nolimits_{\cal L} Y\geq
\min\{mN+1,\mu_{l+1}-l+1\}.
\end{equation}

{\bf Proof.} Obviously $\mu_j\geq jm$ for each
$j=1,\dots,l+1$. Thus
$$
(\mu_{j+1}-j)(N-j)+1\geq \varepsilon(j)+mN+1,
$$
where
$\varepsilon(t)=-(m-1)t^2+(Nm-N-m)t$ has two roots, $t=0$
and $t=N-1-\frac{1}{m-1}$. Thus
$\varepsilon(0)=0$ and $\varepsilon(j)\geq 0$ for
$j=1,\dots,N-2$. Therefore we may omit in (\ref{a5})
the values $j=1,\dots,l-1$. If $l\leq N-2$, we may also
omit the value $j=l$. Q.E.D. for the corollary.


\subsection{Start of the proof of Proposition 2.1}

Let $U_{sm}\subset{\cal H}$ be a non-empty Zariski open
subset, consisting of all collections $(f_*;g_*)$ such that:

(i) $Q=Q(f_*)\subset{\mathbb P}$ is a smooth complete intersection;

(ii) all the divisors $W_i|_Q$ are smooth and the divisor
$$
(W_1+\dots+W_m)|_Q
$$
has normal crossings. In particular, the points of class
$e\geq 1$ make a smooth quasi-projective variety of
codimension $e$ (its closure is the set of points of class
$\geq e$), and the very variety $V=V(f_*;g_*)$ is smooth.

Let $e\in\{ 0,\dots,m\}$ be fixed. Consider the closed
subset
$$
Y_e=\{(x,(f_*;g_*))\in {\mathbb P}\times U_{sm}\,|\,
x\in Q(f_*)\,\,\,\mbox{is non-regular of class}\,\,\, e\}.
$$
Let $\pi\colon{\mathbb P}\times U_{sm}\to U_{sm}$ be
the projection onto the second factor. Now Proposition 2.1
follows from 
\vspace{0.3cm}

{\bf Proposition 2.3.} {\it The closure
$$
\overline{\pi(Y_e)}\subset U_{sm}
$$
is a proper closed subset for any $e\in\{0,1,\dots,m\}$.
}\vspace{0.3cm}

{\bf Proof.} We show that  $\overline{\pi(Y_e)}$ has a
positive codimension in $U_{sm}$. Set
$$
\begin{array}{l} 
\displaystyle
Y_e(x)=Y_e\cap (\{x\}\times U_{\rm sm})\subset
{\mathbb P}\times U_{\rm sm}, \\   \\
\displaystyle
I=\{(x,(f_*))|x\in Q(f_*)\}\subset
{\mathbb P}\times U_{\rm sm}, \\    \\
\displaystyle
I_e=\{(x,(f_*))|x\in Q(f_*)\,\,\mbox{is of class e}\}\subset
{\mathbb P}\times U_{\rm sm}, \\    \\
\displaystyle
I_e(x)=I_e\cap (\{x\}\times U_{\rm sm})\subset
{\mathbb P}\times U_{\rm sm}.
\end{array}
$$
Identifying $\{x\}\times U_{sm}\cong U_{sm}$,
we think of $Y_e(x)$ and $I_e(x)$ as subsets in $U_{sm}$.

In the non-homogeneous presentation with respect to the
system of affine
coordinates $(z_1,\dots,z_{M+k})$ the collection
$(f_*;g_*)$ can be identified with the set of
{\it homogeneous} polynomials $q_{i,j}$, $w_{i,j}$. Among
the polynomials (\ref{a3},\ref{a4}) that appear in the regularity condition
there are precisely $k+e$ linear forms: these are
$$
q_{1,1},\dots,q_{k,1},w_{1,1},\dots,w_{e,1}.
$$
Since $(f_*;g_*)\in U_{sm}$, these linear forms
are linearly independent. Set
$$
P=P(f_*;g_*)=\{v\in{\mathbb C}^{M+k}\,|\,
q_{1,1}=\dots=w_{e,1}=0\},
$$
$P\cong {\mathbb C}^{M-e}$. If the point $x$ is non-regular,
then the set of common zeros of the rest of polynomials in the
set (\ref{a3}) or (\ref{a4}) is of a ``wrong'' codimension. It is easy to
compute that there is not more than $M-e-1$  polynomials of
degree $\geq 2$ in the list (\ref{a3}) or (\ref{a4}). Let us apply Corollary 2.1 to them.


\subsection{The polynomials $h_{i,j}$ depend on each other} 

Note that the polynomials $h_{i,j}$ formally depend on each other, so that generally speaking
we cannot apply Corollary 2.1. However, the
explicit form of the polynomials $h_{i,j}$ (\ref{a2}) makes it
possible to circumwent this obstruction. Indeed, let us
assume that the polynomials
$$
w_{i,1},\dots,w_{i,l_i}
$$ 
for $i\geq e=1$ are fixed.

Let us construct by induction a sequence of homogeneous
polynomials
$$
\xi_{i,l_i+1},\dots,\xi_{i,2l_i},
$$
setting
$$
\begin{array}{l}
\displaystyle
\xi_{i,l_i+1}=-A_{l_i}(w_{i,1},\dots,w_{i,l_i}),\\  \\ \displaystyle
\xi_{i,j+1}=-A_{j}(w_{i,1},\dots,w_{i,l_i},
\xi_{i,l_i+1},\dots,\xi_{i,j}),
\end{array}
$$
$j=l_i+1,\dots,2l_i-1$. It is obvious that the polynomials
$\xi_{i,j}$ depend on $w_{i,1},\dots,w_{i,l_i}$ only and
therefore are also fixed. 
\vspace{0.1cm}

{\bf Lemma 2.3.} {\it For any closed subset $T\subset P$ the
following conditions are equivalent:

\noindent
{\rm (i)} all the polynomials $h_{i,j}$ vanish on $T$,
$j=l_i+1,\dots, a\leq 2l_i$;

\noindent
{\rm (ii)} all the polynomials $w_{i,j}-\xi_{i,j}$ vanish
on $T$, $j=l_i+1,\dots, a\leq 2l_i$.
}
\vspace{0.1cm}

{\bf Proof.} Indeed,
$$
h_{i,l_i+1}\equiv w_{i,l_i+1} -\xi_{i,l_i+1},
$$
so that for $a=l_i+1$ the claim of the lemma is obvious.
Now we argue by induction on $a\geq l_i+2$. If the lemma
is true for $a\leq c\leq 2l_i-1$, and any of the
conditions (i), (ii) holds for $a=c+1$, then in any case
$$
(w_{i,j}-\xi_{i,j})|_T\equiv 0
$$
for $j=l_i+1,\dots,c$. Therefore
$$
w_{i,j}|_T\equiv \xi_{i,j})|_T
$$
for $j=l_i+1,\dots,c$, so that the conditions
$$
[w_{i,c+1}+A_c(w_{i,1},\dots,w_{i,l_i},
w_{i,l_i+1},\dots,w_{i,c})]|_T\equiv 0
$$
and 
$$
\begin{array}{rcccl}
\displaystyle
[   &   w_{i,c+1}  &  +  &  A_c(w_{i,1},\dots,w_{i,l_i},
w_{i,l_i+1},\dots,w_{i,c})  &   ]|_T\equiv 0 \\  \\
\displaystyle
  &     &    &  \|  &    \\   \\ \displaystyle
  &     &    &   -\xi_{i,c+1}  &     
\end{array}
$$
are equivalent. Q.E.D. for the lemma.

Thus for any fixed set of polynomials
$w_{i,1},\dots,w_{i,l_i}$ the polynomials $h_{i,j}$
in the regularity condition can be replaced by the
polynomials $w_{i,j}-\xi_{i,j}$. However, the polynomials
$w_{i,j}$, $j=l_i+1,\dots,2l_i$, are arbitrary and
therefore the polynomials $(w_{i,j}-\xi_{i,j})$ are
also arbitrary: essentially we just shift the origin
in the space of homogeneous polynomials of degree $j$
in $z_*$. Thus when we estimate the codimension 
$\mathop{\rm codim}\nolimits_{I(x)}Y(x)$ that comes
from the regularity condition being not satisfied we
may apply Corollary 2.1 as if all the polynomials $q_{i,j}$,
$h_{i,j}$ were homogeneous polynomials independent of each other. Now let us, at long last, estimate this
codimension.


\subsection{Estimate for the codimension}

Set $Y_1(x)=Y^+(x)\cup Y^{\sharp}(x)$, where
$(f_*;g_*)\in Y^+(x)$ if and only if $l_1\geq 3$,
whereas $(f_*;g_*)\in Y^{\sharp}(x)$  if and only if
$l_1=2$.

Obviously,
\begin{equation}
\label{a11}
\pi(Y_e)=\mathop{\bigcup}\limits_{x\in{\mathbb P}}Y_e(x).
\end{equation}
It suffices to show that
\begin{equation}
\label{a12}
\mathop{\rm codim}\nolimits_{U_{sm}}\pi(Y_e)\geq 1.
\end{equation}
By (\ref{a11}) this estimate follows immediately from
the inequality
$$
\mathop{\rm codim}\nolimits_{U_{sm}}Y_e(x)\geq M+k+1,
$$
which is what we shall actually prove. Let us consider
each possible case in turn.
\vspace{0.3cm}

{\bf Case $e\geq 2$.} Here by Definition 2.2 (ii) we have
$$
\sum^k_{i=1}(d_i-1)+\sum^m_{i=e+1}l_i=
M-\sum^e_{i=1}l_i\leq M-e-2
$$
homogeneous polynomials of degree $\geq 2$ on 
$P\cong{\mathbb C}^{M-e}$. According to Corollary 2.1,
$$
\mathop{\rm codim}\nolimits_{I_e(x)}Y_e(x)\geq 2M-2e-1.
$$
Since obviously
$$
\mathop{\rm codim}\nolimits_{U_{sm}}I_e(x)=k+e,
$$
we get finally the estimate
$$
\mathop{\rm codim}\nolimits_{U_{sm}}Y_e(x)\geq 2M+k-e-1.
$$
Taking the union over all $x\in{\mathbb P}$, we get
$$
\begin{array}{ccc}
\displaystyle
   &  \overline{
\mathop{\bigcup}\limits_{x\in{\mathbb P}}Y_e(x)}  &  \\ \\
   &    \displaystyle  \|  &   \\  \\
\displaystyle
\mathop{\rm codim}\nolimits_{U_{sm}} &  \overline{Y_e} & \geq M-e-1.
\end{array}
$$
Since $e<M/2$, the inequality (\ref{a12}) is proved, so that
$\overline{Y_e}\subset U_{sm}$ is a proper 
closed subset.
\vspace{0.3cm}

{\bf Case $e=1$.} In the $+$-subcase we proceed as
above for $e\geq 2$. Let consider the ${\sharp}$-subcase.
Here we have precisely $M-2$ polynomials on 
${\mathbb C}^{M-1}$, so that
\begin{equation}
\label{a13}
\mathop{\rm codim}\nolimits_{I(x)}Y^{\sharp}(x)\geq
\min(2M-3,\beta),
\end{equation}
where
$$
\begin{array}{ccc}
\displaystyle
\beta  &  =   &  \displaystyle
\sum^k_{i=1}\sum^{d_i}_{j=2}(j-1)+
\sum^m_{i=2}\sum^{2l_i}_{j=l_i+1}(j-1)+2 \\  \\
   &    &   \|   \\   \\
\displaystyle
   &    &  \displaystyle
1/2\left[ 
\sum^k_{i=1}d_i(d_i-1)+
\sum^m_{i=2}l_i(3l_i-1)
\right] + 2.
\end{array}
$$
If the minimum in (\ref{a13}) is attained at $2M-3$,
then we proceed as in the $+$-case above. Thus it suffices
to consider the case when the minimum is attained at
$\beta$. Let us estimate $\beta$ as a function of
non-negative real varieties $d_i$, $l_i$, satisfying
the constraints
$$
\sum^k_{i=1}d_i+\sum^m_{i=2}l_i=M+k-2, \quad
d_i\geq 2,\quad l_i\geq 2.
$$
\vspace{0.1cm}

{\bf Lemma 2.4.} {\it The following inequality holds:}
$$
\beta\geq M.
$$
\vspace{0.1cm}

First of all, assuming the claim of the lemma,
let us complete the ${\sharp}$-case.
We have
$$
\mathop{\rm codim}\nolimits_{I_1(x)}Y^{\sharp}(x)\geq M,\quad
\mathop{\rm codim}\nolimits_{U_{sm}}I_1(x)=k+1,
$$
so that
$$
\mathop{\rm codim}\nolimits_{U_{sm}}Y^{\sharp}(x)\geq M+k+1,
$$
and taking the union over all the points $x\in{\mathbb P}$
we get
$$
\mathop{\rm codim}\nolimits_{U_{sm}}\overline{Y^{\sharp}}\geq 1,
$$
which is what we need.
\vspace{0.1cm}

{\bf Proof of Lemma 2.4.} For convenience of computations
we start with the following auxiliary claim.
\vspace{0.1cm}

{\bf Lemma 2.5.} (i) {\it For $s_i\geq 2$, 
$\sum\limits^c_{i=1}s_i=B\geq 2c$, where $c\in{\mathbb Z}_+$ is a
fixed positive integer the following inequality
holds:}
$$
\sum^c_{i=1}s_i(3s_i-1)\geq 5B.
$$

(ii) {\it For $s_i\geq 2$, 
$\sum\limits^k_{i=1}s_i=A\geq 2k$, where $k\in{\mathbb Z}_+$ is a
fixed positive integer the following inequality
holds:}
$$
\sum^k_{i=1}s_i(s_i-1)\geq A(\frac{\displaystyle A}{
\displaystyle k}-1).
$$

{\bf Proof}: elementary computations. It is easy to see
that in both cases the minimum is attained at
$s_1=\dots=s_c$ (or $=s_k$.)

Setting
$$
\sum^k_{i=1} d_i=A,\quad
\sum^m_{i=2} l_i=B,
$$
we obtain by Lemma 2.5 the estimate
$$
\beta\geq\frac12 (A(\frac{\displaystyle A}{
\displaystyle k}-1)+5B)+2.
$$
Now to prove Lemma 2.4 it is sufficient to check that
the inequality
$$
\varepsilon(A,B)=A(\frac{\displaystyle A}{
\displaystyle k}-1)+5B+4\geq 2M
$$
is true under the constraints
$$
1\leq k\leq \frac{M-1}{2},\quad A\geq 2k,\quad
B\geq 0,\quad A+B=M+k-2.
$$
Replacing $B$ by $M+k-2-A$, let us consider the function
of real variable $A\in{\mathbb R}_+$
$$
\zeta(A)=A(\frac{\displaystyle A}{
\displaystyle k}-6).
$$
Its minimum on the interval $I=[2k,M+k-2]$ is attained
either at $A=3k$ (if $2k\leq M-2$), or at
$A=M+k-2$ (if $2k=M-1$). In the first case we get
$$
\varepsilon(A,B)\geq 2A+5B+6=2(M+k-2)+3B+6\geq 2M+4.
$$
In the second case $B=0$ and by elementary computations
we get
$$
\varepsilon(A,B)\geq (M+k-2)\frac{
\displaystyle 2k-1}{\displaystyle
k}+6\geq 2M.
$$
This completes the proof of Lemma 2.4.
\vspace{0.3cm}

{\bf Case $e=0$.} Here to prove the estimate (\ref{a12}) it is
sufficient to show that
\begin{equation}
\label{a14}
\beta=\frac12\left[
\sum^k_{i=1}d_i(d_i-1)+
\sum^{m-1}_{i=1}l_i(3l_i-1)+l_m(3l_m-5)
\right]+2\geq M+1,
\end{equation}
since if this is the case then the estimate
$$
\mathop{\rm codim}\nolimits_{U_{sm}}Y_0(x)\geq M+k+1
$$
holds, and arguing as above we see that 
$\overline{Y_0}\subset U_{sm}$ is a proper
closed subset.

Consider $\beta$ as a function of non-negative real
variables $d_i$, $l_i$, and set
$$
\sum^k_{i=1} d_i=A,\quad
\sum^{m-1}_{i=1} l_i=B.
$$
Assuming $A,B$ to be fixed, we get by Lemma 2.5
the estimate
$$
\beta\geq\frac12 (A(\frac{\displaystyle A}{
\displaystyle k}-1)+5B+l_m(3l_m-5))+2.
$$
Thus to prove the estimate (\ref{a14}) it is sufficient to
check that the inequality
\begin{equation}
\label{a15}
A(\frac{\displaystyle A}{
\displaystyle k}-1)+5B+l_m(3l_m-5)+2\geq 2M.
\end{equation}
First of all let us get rid of $l_m$. Since
$$
A+B+l_m=M+k, \quad l_m\geq 2,
$$
we assume $A$ and $B+l_m$ to be fixed. Since the
derivative
$$
[-5t+t(3t-5)]'=6t-10
$$
is positive for $t\geq 2$, the minimum of the 
left-hand side  of (\ref{a15}) is attained at
$l_m=2$. Therefore it is sufficient to prove the
inequality
$$
A(\frac{\displaystyle A}{
\displaystyle k}-1)+5B+4\geq 2M
$$
for $A+B=M+k-2$ and the standard constraints for
$A,B$ and $k$. But this has already been done when
the ${\sharp}$-case was considered.

Proof of Proposition 2.1 is complete. For a general
collection $(f_*;g_*)\in{\cal H}$ each point
$x\in Q(f_*)$ is regular.


\section{Hypertangent divisors}

\subsection{How to obtain a bound for the multiplicity}

Let $X$ be a smooth projective variety, $H\in\mathop{\rm Pic} X$ an
ample class. For an irreducible subvariety $Y\subset X$
its $H$-degree (or, simply, degree, when it is clear
what ample class is meant) is the integer
$$
\mathop{\rm deg}\nolimits_H Y=(Y\cdot H^{\mathop{\rm dim} Y}).
$$
By linearity the $H$-degree is defined for any cycle
(which is assumed to be equidimensional).

The symbol
$$
\frac{\displaystyle \mathop{\rm mult}\nolimits_x}{\displaystyle \mathop{\rm deg}\nolimits_H} Y
$$
means the ratio $(\mathop{\rm mult}\nolimits_x Y)/\mathop{\rm deg}\nolimits_H Y$, where $x\in X$ is
a point. Set
$$
\lambda_e(x)=
\mathop{\rm sup}\limits_{
\begin{array}{c} \scriptstyle
T\subset X, \\ \scriptstyle
\mathop{\rm codim}\nolimits_X T=e
\end{array}} \left\{
\frac{\displaystyle \mathop{\rm mult}\nolimits_x}{\displaystyle \mathop{\rm deg}\nolimits_H} T\right\},
$$
where the supremum is taken over all irreducible
subvarieties of codimension $e\geq 1$. 

Assume that on $X$ there exists a set of effective
divisors 
$$
D_i\in |a_i H|,
$$
$i=1,\dots,N$, such that the set-theoretic intersection
$$
D_1\cap \dots\cap D_N
$$
is of codimension $N\leq \mathop{\rm dim} X$ in a neighborhood of the
point $x$. Set
$$
\mu_i=\mathop{\rm mult}\nolimits_x D_i\geq 1.
$$ 
Let $T\subset X$ be an irreducible subvariety of
codimension $e\geq 1$, where $N\geq e+1$ and $T\ni x$.
\vspace{0.1cm}

{\bf Lemma 3.1.}{\it There exists a subset
${\cal L}\subset\{1,\dots,N\}$ of cardinality $N-e$
(after re-numbering we assume, to simplify the
notations, that ${\cal L}=\{1,\dots,N-e\}$) and a
sequence of irreducible subvarieties $T_i$,
$i=0,1,\dots,N-e$, such that:

{\rm (i)} $\mathop{\rm codim}\nolimits T_i=e+i$;

{\rm (ii)} $T_0=T$, $T_i\not\subset D_i$ and 
$T_{i+1}$ is an irreducible component of the 
effective cycle $T_i\cap D_i$;

{\rm (iii)} $T_i\ni x$ and the following inequality
holds:
\begin{equation}
\label{b1}
\frac{\displaystyle \mathop{\rm mult}\nolimits_x}{\displaystyle \mathop{\rm deg}\nolimits_H} T_i\geq
\frac{\displaystyle
\mu_i}{\displaystyle a_i}\cdot
\frac{\displaystyle \mathop{\rm mult}\nolimits_x}{\displaystyle \mathop{\rm deg}\nolimits_H} T_{i-1}
\end{equation}
for all $i=1,\dots,N-e$.
}\vspace{0.1cm}

{\bf Proof.} Let us prove the existence of ${\cal L}$
and the set of subvarieties $T_i$ by induction on
$i\in\{0,1,\dots,N-e\}$. For $i=0$ we have nothing
to prove.
\vspace{0.1cm}

{\bf Lemma 3.2.} {\it There is a divisor $D_i$, $1\leq i\leq N$,
such that} $T\not\subset D_i$.
\vspace{0.1cm}

{\bf Proof.} Assume the converse: $T\subset D_i$ for all
$i=1,\dots,N$. Then
$$
T\subset D_1\cap\dots\cap D_N
$$
so that
$$
\mathop{\rm codim}\nolimits_x(D_1\cap\dots\cap D_N)\leq \mathop{\rm codim}\nolimits T=e\leq N-1
$$
contrary to our assumption about the set $D_1,\dots, D_N$.
Q.E.D. for the lemma.

To simplify the notations we assume that $T\not\subset D_1$.
Obviously,
$$
\begin{array}{rcl}
\displaystyle
\mathop{\rm mult}\nolimits_x(T\circ D_1)  &  \geq  &  
\mathop{\rm mult}\nolimits_x T\cdot \mathop{\rm mult}\nolimits_x D_1, \\  \\
\displaystyle
\mathop{\rm deg}\nolimits_H(T\circ D_1)  &  =  &  a_1 \mathop{\rm deg}\nolimits_H T.
\end{array}
$$
Therefore, the following estimate holds:
$$
\frac{\displaystyle \mathop{\rm mult}\nolimits_x}{\displaystyle \mathop{\rm deg}\nolimits_H} (T\circ D_1)\geq 
\frac{\displaystyle
\mu_1}{\displaystyle a_1}\cdot
\frac{\displaystyle \mathop{\rm mult}\nolimits_x}{\displaystyle \mathop{\rm deg}\nolimits_H} T.
$$
The inequality
\begin{equation}
\label{b2}
\frac{\displaystyle \mathop{\rm mult}\nolimits_x}{\displaystyle \mathop{\rm deg}\nolimits_H} Y\geq \gamma
\end{equation}
is certainly non-linear in $Y$. However, it is
equivalent to the linear inequality
$$
\mathop{\rm mult}\nolimits_x Y\geq \gamma \mathop{\rm deg}\nolimits_H Y.
$$
Therefore if (\ref{b2}) holds for an effective cycle
$Y$, then there exists a component $Y^+$ of this cycle,
that is, an irreducible
subvariety in $X$, such that
$$
\frac{\displaystyle \mathop{\rm mult}\nolimits_x}{\displaystyle \mathop{\rm deg}\nolimits_H} Y^+\geq \gamma
$$
(since the $H$-degree of an irreducible subvariety
is always strictly positive). So we get that there is
an irreducible component $T_1$ of the effective cycle
$(T\circ D_1)$, an irreducible subvariety of codimension
$e+1$, such that
$$
\frac{\displaystyle \mathop{\rm mult}\nolimits_x}{\displaystyle \mathop{\rm deg}\nolimits_H} T_1\geq
\frac{\displaystyle
\mu_1}{\displaystyle a_1}\cdot
\frac{\displaystyle \mathop{\rm mult}\nolimits_x}{\displaystyle \mathop{\rm deg}\nolimits_H} T,
$$
which is what we need for $i=1$.

Assume that the subset $\{1,\dots,j\}$, $j\leq N-e-1$,
and a sequence of irreducible subvarieties $T_1,\dots,T_j$
satisfy the conditions (i)-(iii).
\vspace{0.1cm}

{\bf Lemma 3.3.} {\it There is a divisor $D_i$, 
$j+1\leq j\leq N$, such that} $T_j\not\subset D_i$.
\vspace{0.1cm}

{\bf Proof.} Assume the converse: $T_j\subset D_i$ for
all $i=j+1,\dots,N$. Then
$$
T_j\subset  D_{j+1}\cap \dots\cap D_N.
$$
Taking into consideration that by construction
$$
T_j\subset  D_{1}\cap \dots\cap D_j,
$$
we get
$$
T_j\subset  D_{1}\cap \dots\cap D_N.
$$
However $x\in T_j$ and the codimension of the
subvariety $T_j$ is equal to
$$
e+j\leq N-1,
$$
which gives again (as in the proof of Lemma 3.2) a
contradiction with what we assumed about the
collection $D_1,\dots,D_N$. Q.E.D. for the lemma.

After re-numbering we may assume that 
$T_j\not\subset D_{j+1}$. Now we argue as above:
$$
\frac{\displaystyle \mathop{\rm mult}\nolimits_x}{\displaystyle \mathop{\rm deg}\nolimits_H} (T_j\circ D_{j+1}) \geq
\frac{\displaystyle
\mu_{j+1}}{\displaystyle a_{j+1}}\cdot
\frac{\displaystyle \mathop{\rm mult}\nolimits_x}{\displaystyle \mathop{\rm deg}\nolimits_H} T_j
$$
and therefore there is an irreducible component $T_{j+1}$
of the effective cycle $(T_j\circ D_{j+1})$ that satisfies
the inequality (\ref{b1}). Proof of Lemma 3.3 is complete.
\vspace{0.3cm}

{\bf Corollary 3.1.} {\it The following inequality holds}
\begin{equation}
\label{b3}
\lambda_e(x)\cdot
\mathop{\rm min}\limits_{
\begin{array}{c}
\scriptstyle
{\cal L}\subset\{1,\dots,N\} \\
\scriptstyle
\sharp {\cal L}=N-e
\end{array}
}\left(
\prod_{i\in{\cal L}}\frac{\displaystyle
\mu_i}{\displaystyle a_i}
\right)\leq \lambda_N(x).
\end{equation}
\vspace{0.3cm}

{\bf Proof.} In the notations above
$$
\lambda_N(x)\geq \frac{\displaystyle \mathop{\rm mult}\nolimits_x}{\displaystyle \mathop{\rm deg}\nolimits_H} T_N\geq
\left(
\prod_{i\in{\cal L}}\frac{\displaystyle
\mu_i}{\displaystyle a_i}
\right)\frac{\displaystyle \mathop{\rm mult}\nolimits_x}{\displaystyle \mathop{\rm deg}\nolimits_H} T.
$$
Here ${\cal L}\subset\{1,\dots,N\}$ is a subset
of cardinality $N-e$, which depends, generally
speaking, on $T$. The more so,
$$
\lambda_N(x)\geq
\mathop{\rm min}\limits_{
\begin{array}{c}
\scriptstyle
{\cal L}\subset\{1,\dots,N\} \\
\scriptstyle
\sharp {\cal L}=N-e
\end{array}
}\left(
\prod_{i\in{\cal L}}\frac{\displaystyle
\mu_i}{\displaystyle a_i}
\right)\cdot \frac{\displaystyle \mathop{\rm mult}\nolimits_x}{\displaystyle \mathop{\rm deg}\nolimits_H} T.
$$
The first factor in the right-hand side does not
depend on $T$. Since the variety $T\subset X$ is
absolutely arbitrary, we get the inequality
(\ref{b3}). Q.E.D. for the corollary.
\vspace{0.3cm}

{\bf Corollary 3.2.} {\it  Assume that the linear system
$|H|$ is free and defines a (finite) morphism
$\varphi_{|H|}\colon X\to{\mathbb P}^k$. Then the following
estimate holds:}
\begin{equation}
\label{b4}
\lambda_e(x)\leq
\left(
\mathop{\rm min}\limits_{
\begin{array}{c}
\scriptstyle
{\cal L}\subset\{1,\dots,N\} \\
\scriptstyle
\sharp {\cal L}=N-e
\end{array}
}\left(
\prod_{i\in{\cal L}}\frac{\displaystyle
\mu_i}{\displaystyle a_i}
\right)
\right)^{-1}.
\end{equation}
\vspace{0.3cm}

{\bf Proof.} For any irreducible subvariety $T\subset X$
of codimension $e\geq 1$ there exist divisors
$D_i\in |H|$, $i=1,\dots,\mathop{\rm dim} T$, such that:

\begin{itemize}
\item
$D_i\ni x$, in particular 
$\mathop{\rm mult}\nolimits_x D_i\geq 1$;

\item
the intersection 
$$
T^{\sharp}=T\cap D_1\cap \dots\cap D_{\mathop{\rm dim} T}
$$
is zero-dimensional. 

\end{itemize}

\noindent
Obviously,
$$
\mathop{\rm deg}\nolimits T^{\sharp}=\mathop{\rm deg}\nolimits T
$$
and 
$$
\mathop{\rm mult}\nolimits_x T^{\sharp}\geq
\mathop{\rm mult}\nolimits_x T\cdot
\prod^{\mathop{\rm dim} T}_{i=1}\mathop{\rm mult}\nolimits_xD_i\geq
\mathop{\rm mult}\nolimits_x T.
$$
However, $T^{\sharp}$ is a zero-dimensional scheme,
so that
$$
\mathop{\rm deg}\nolimits T^{\sharp}\geq \mathop{\rm mult}\nolimits_x T^{\sharp}.
$$
From this inequality we obtain that $\lambda_e(x)\leq 1$
for all $x$ and $e$. Now applying the previous corollary
we complete the proof.


\subsection{Construction of hypertangent divisors}

To realize the method described in Sec. 3.1 for the
iterated double covers, let us first of all fix some notations. As above, we have a system of affine
coordinates $(z_1,\dots,z_{M+k})$ with the origin at
the point $p\in Q$. For $1\leq i\leq k$, 
$1\leq j\leq d_i-1$ set
$$
f_{i,j}=q_{i,1}+\dots+q_{i,j},
$$
$$
D^{{\mathbb P}}_{i,j}=\overline{\{f_{i,j}=0\}}
$$
(the closure is taken in ${\mathbb P}$),
$$
D^{Q}_{i,j}=D^{{\mathbb P}}_{i,j}|_Q,
$$
$$
D^{f}_{i,j}=\sigma^{-1}(D^{Q}_{i,j}).
$$
Assume that the point $p\in Q$ is of class $e\geq 0$,
and, moreover, if $e\geq 1$, then 
$p\in W_1\cap\dots\cap W_e$.
Set
$$
D^{+}_{i}=\overline{\{w_{i,1}=0\}}|_Q,
$$
$$
D_i=\sigma^{-1}(D^+_i),
$$
$i\in\{1,\dots,e\}$. Finally, for a point $x\in V$ such
that $\sigma(x)=p$, let us define with respect to the
coordinates $u_i$, $i\geq e+1$ (see Sec. 2.1), the
following divisors:
$$
D^{g}_{i,j}=\overline{\{u_i-[\sqrt{g_i}]_j=0\}}|_V,
$$
$$
D^{+}_{i,j}=\sigma(D_{i,j}),
$$
where $j=l_i,\dots,2l_i-1$. It is easy to see that
$$
D^{?}_{i,j}\in |jH|,\quad D_i\in |H|
$$
for any $i,j$, listed above, where $?\in \{f,g\}$.
\vspace{0.1cm}

{\bf Lemma 3.4.} {\it  {\rm (i)} For any $i,j$ the following
inequality holds:
\begin{equation}
\label{b5}
\mathop{\rm mult}\nolimits_x D_{i,j}\geq j+1.
\end{equation}

\noindent
{\rm (ii)} For any $i\in\{1,\dots,e\}$, where $e\geq 1$,
the following inequality holds:
$$
\mathop{\rm mult}\nolimits_x D_i\geq 2.
$$
}
\vspace{0.1cm}

{\bf Proof.} (i) To begin with, take $1\leq i\leq k$.
Obviously,
$$
f_{i,j}|_Q=(-q_{i,j+1}-\dots-q_{i,d_i})|_Q,
$$
since $f_i|_Q\equiv 0$, which implies the estimate
(\ref{b5}). Now assume that $i\geq e+1$. By the definition
of the class of a point, $g_i(p)\neq 0$. In the open
affine subset $U\subset{\mathbb P}^{\sharp}$, 
$U\cong{\mathbb C}^{M+k+m}$ with coordinates $(z_*,u_*)$
we get
$$
[u_i^2-g_i(1,z_1,\dots,z_{M+k})]|_V\equiv 0.
$$
Since obviously
$$
(y_i+[\sqrt{g_i}]_j)(x)\neq 0,
$$
we obtain that locally the divisor $D_{i,j}$ is given by
the equation
$$
g^{(j)}_i|_V=0.
$$
As we have seen above, the first non-zero homogeneous component in $g^{(j)}_i$ is $h_{j+1}[g_i]$; it is of
degree $j+1$.

This completes the proof of the first part of the lemma.

(ii) Assume that $e\geq 1$. The hyperplane $\{w_{i,1}=0\}$
is tangent to the hypersurface $W_i$ at the point $p$.
Hence the divisor $D_i$ is singular at the point $x$. This
is what we need.


\subsection{The regularity condition for hypertangent divisors}

To apply the techniques of Sec. 3.1 it is not enough just
to know the multiplicities of hypertangent divisors at the
point $x$. We need more precise information about the
tangent cones to these divisors. Set
$$
E=T_p(Q\cap W_1\cap\dots\cap W_e), 
\quad
U=\langle u_1,\dots,u_e\rangle^*.
$$
As we have seen above (Lemma 2.1), the following isomorphism
holds
\begin{equation}
\label{b6}
T_q V\cong E\oplus U.
\end{equation}
The subspace $E\subset T_p{\mathbb P}={\mathbb C}^{M+k}_{(z_*)}$
is given by the system of linear equations
$$
q_{1,1}=\dots=q_{k,1}=w_{1,1}=\dots=w_{e,1}=0.
$$
The local computations performed in the proof of the
previous lemma show that by the isomorphism (\ref{b6})
the tangent cones to the hypertangent divisors are
given by the following equations:

\begin{itemize}
\item
to the divisors $D^{f}_{i,j}$ ---
\begin{equation}
\label{b7}
q_{i,j+1}|_E=0
\end{equation}
for $1\leq i\leq k$;

\item
to the divisors $D^{g}_{i,j}$ ---
\begin{equation}
\label{b8}
h_{j+1}[g_i]|_E=0
\end{equation}
for $i\geq e+1$;

\item
to the divisors $D_i$ ---
\begin{equation}
\label{b9}
u^2_i=w_{i,2}|_E
\end{equation}
for $1\leq i\leq e$, if $e\geq 1$.

\end{itemize}

Indeed, the equations (\ref{b7}) have been obtained
above in Sec. 3.2. The equations (\ref{b8}) and
(\ref{b9}) follow from the local computations made in
Sec. 3.2 if one takes into consideration that any 
linear form $L(z_*)$ that vanish on $E$ defines an
element in the {\it square} of the maximal ideal
${\cal M}_{x,V}$ of the point $x$ in the local ring
${\cal O}_{x,V}$:
$$
\sigma^*(L(z_*)|_Q)\in {\cal M}^2_{x,V}.
$$
Thus if a pair of homogeneous polynomials
$P^+(z_*)$, $P^-(z_*)$ of degree $a\geq 1$ coincide on
$E$, that is, $(P^+-P^-)|_E\equiv 0$, then
$$
\sigma^*(P^+(z_*)|_Q)\equiv \sigma^*(P^-(z_*)|_Q)
\mathop{\rm mod} {\cal M}^{a+1}_{x,V}.
$$
From here the equations (\ref{b8}) and
(\ref{b9}) follow immediately.
\vspace{0.1cm}

{\bf Lemma 3.5.} {\it The set ${\cal D}=
\{D^{f}_{i,j},D_i,D^{g}_{i,j}\}$ of all hypertangent
divisors satisfies the regularity condition}
$$
\mathop{\rm codim}\nolimits_x\mathop{\bigcap}\limits_{D\in{\cal D}} D=
\sharp {\cal D}.
$$
\vspace{0.1cm}

{\bf Proof.} It is sufficient to compute the
codimension
\begin{equation}
\label{b10}
\mathop{\rm codim}\nolimits_{T_xV} \mathop{\bigcap}\limits_{D\in{\cal D}} T_xD.
\end{equation}
Since $T_xV=U\oplus E$ and the coordinates $u_i$ come
into the equations (\ref{b9}) only --- each in its own,
the codimension (\ref{b10}) is equal to
\begin{equation}
\label{b11}
\mathop{\rm codim}\nolimits_E\{
q_{i,j+1}|_E=0,\, h_{j+1}[g_i]=0\},
\end{equation}
where the indices $i,j$ comprise the sets
$$
\{1\leq i\leq k,\, 1\leq j\leq d_i-1\}
$$
and 
$$
\{e+1\leq i\leq m,\, l_i+1\leq j\leq 2l_i\},
$$
respectively. Adding the equations of the hyperplane $E$,
we see that the codimension (\ref{b10}) is precisely the
codimension in ${\mathbb C}^{M+k}$ of the set determined by
all the polynomials that come into the regularity condition
(\ref{a3}) or (\ref{a4}). The claim of Lemma 3.5 follows 
from this fact immediately.


\subsection{The Lefschetz theorem once again}

Assume that the point $p=\sigma(x)$ is of class $e=0$.
\vspace{0.1cm}

{\bf Lemma 3.6.} {\it The set-theoretic intersection
$$
T=D_{1,1}\cap\dots\cap D_{k,1}\subset V
$$
is of codimension $k$, coincides with the scheme-theoretic
intersection
$$
T=(D_{1,1}\circ\dots\circ D_{k,1})
$$
and satisfies the equalities}
$$
\mathop{\rm deg}\nolimits T=\mathop{\rm deg}\nolimits V, \quad \mathop{\rm mult}\nolimits_x T=2^k.
$$
\vspace{0.1cm}

{\bf Proof.} Let show by induction on $i=1,\dots,k$, that
the set-theoretic intersection
$$
T_i=\mathop{\bigcap}\limits_{j=1}^{i} D_{j,1}\subset V
$$
is of codimension $i$, coincides with the scheme-theoretic
intersection:
$$
T_i=(D_{1,1}\circ\dots\circ D_{i,1})
$$
and satisfies the equalities
$$
\mathop{\rm deg}\nolimits T_i=\mathop{\rm deg}\nolimits V, \quad \mathop{\rm mult}\nolimits_x T_i=2^i.
$$
Indeed, for $i=1$ it is true in an obvious way: the
tangent cone
$$
T_x D_{1,1}\subset T_x V\cong T_p Q
$$
is given by the quadratic equation
$$
q_{1,2}|_{\{q_{1,1}=\dots=q_{k,1}=0\}}=0,
$$
which is non-trivial by the regularity condition.
To come over from $i$ to $i+1$, let us use the Lefschetz
theorem: by the regularity condition the set of common
zeros of the system of equations
\begin{equation}
\label{b12}
q_{1,1}=\dots=q_{k,1}=q_{1,2}=\dots=q_{i,2}=
q_{i+1,2}=0
\end{equation}
is of codimension precisely $i+1$ in $T_x V$. Hence
$$
T_x T_i\not\subset  T_x D_{i+1,1}
$$ 
and thus
$$
T_i\not\subset D_{i+1,1}.
$$
But $T_i\subset V$ is an irreducible subvariety of
codimension $i$. Therefore the scheme-theoretic
intersection
$$
T^+_{i+1}=(T_i\circ D_{i+1,1})
$$
is an effective cycle of codimension $i+1$. However,
$$
\mathop{\rm deg}\nolimits T_i=\mathop{\rm deg}\nolimits T^+_{i+1}=\mathop{\rm deg}\nolimits V
$$
and $i+1<\mathop{\rm dim} V/2$, so by the Lefschetz theorem we get
that $T^+_{i+1}=T_{i+1}$ is an irreducible subvariety,
the class of which generates $A^{i+1}V$. Finally, the
system of equations (\ref{b12}) gives precisely the
tangent cone $T_x T_{i+1}$ (as an effective algebraic
cycle, that is, respecting the multiplicities of the
components). Thus we get
$$
\mathop{\rm mult}\nolimits_x T_{i+1}=2^{i+1},
$$
which is what we need.
\vspace{0.1cm}

{\bf Definition 3.1.} An irreducible subvariety 
$x\in Y\subset V$ of codimension 2 is said to be
{\it correct at the point} $x$ (where $p=\sigma(x)$ is
a point of class 0), if there exists an irreducible subvariety $R\subset Y$ of codimension
$$
\mathop{\rm codim}\nolimits_V R=k+1,
$$
satisfying the inequality
\begin{equation}
\label{b13}
\frac{\displaystyle \mathop{\rm mult}\nolimits_x}{\displaystyle \mathop{\rm deg}} R\geq 2^{k-1}\cdot \frac{\displaystyle \mathop{\rm mult}\nolimits_x}{\displaystyle \mathop{\rm deg}} Y.
\end{equation} 
\vspace{0.1cm}

{\bf Lemma 3.7.} {\it If the subvariety $Y\ni x$ is not
correct at the point $x$, then the following estimate
holds:}
$$
\frac{\displaystyle \mathop{\rm mult}\nolimits_x}{\displaystyle \mathop{\rm deg}} Y\leq \frac{\displaystyle 4
}{\displaystyle \mathop{\rm deg}\nolimits V}.
$$
\vspace{0.1cm}

{\bf Proof.} Let us apply Lemma 3.1 to the irreducible
subvariety $Y$ and the set of hyperplane sections
$\{D_{i,1}\}$. We obtain that there exists an
irreducible subvariety $R^{\sharp}\subset Y$ of
codimension $k$ (with respect to $V$) such that
$$
\frac{\displaystyle \mathop{\rm mult}\nolimits_x}{\displaystyle \mathop{\rm deg}} R^{\sharp}\geq 2^{k-2}\cdot \frac{\displaystyle \mathop{\rm mult}\nolimits_x}{\displaystyle \mathop{\rm deg}} Y.
$$
Two cases are now possible. If
$$
R^{\sharp}\neq T=D_{1,1}\cap\dots\cap D_{k,1},
$$
then there is a hyperplane section $D_{a,1}$, 
$1\leq a \leq k$, which does not contain $R^{\sharp}$.
Thus 
$$
\mathop{\rm codim}\nolimits_V (R^{\sharp}\cap D_{a,1})=k+1,
$$
so that $(R^{\sharp}\circ D_{a,1})$ is an effective cycle
of codimension $k+1$, satisfying the inequality
$$
\frac{\displaystyle \mathop{\rm mult}\nolimits_x}{\displaystyle \mathop{\rm deg}} (R^{\sharp}\circ D_{a,1})\geq 2^{k-1}\cdot
\frac{\displaystyle \mathop{\rm mult}\nolimits_x}{\displaystyle \mathop{\rm deg}} Y.
$$
Therefore, there is an irreducible component $R$ of the
effective cycle $(R^{\sharp}\circ D_{a,1})$, satisfying
the inequality (\ref{b13}). Therefore, $Y$ is a correct
subvariety, which contradicts our assumption.

Thus $R^{\sharp}=T=D_{1,1}\cap\dots\cap D_{k,1}$. 
Consequently, by Lemma 3.6 we get the inequality
$$
\frac{\displaystyle \mathop{\rm mult}\nolimits_x}{\displaystyle \mathop{\rm deg}} Y\leq 2^{2-k}\cdot \frac{\displaystyle \mathop{\rm mult}\nolimits_x}{\displaystyle \mathop{\rm deg}} T=
\frac{\displaystyle 4
}{\displaystyle \mathop{\rm deg}\nolimits V},
$$
as we have claimed. Q.E.D. for the lemma.


\subsection{The points of class $e=0$}

Let us prove, at long last, the estimate (\ref{c1}) for a
point $x\in V$ such that its image $p=\sigma(x)\in{\mathbb P}$
is of class $0$. In order to do this, we apply
the construction of Lemma 3.1 to an arbitrary subvariety
$R\subset V$ of codimension $k+1$, $R\ni x$, and to
the set of $M-1$ divisors, which we get putting together
the $f$-collection
$$
\{D^{f}_{i,j}\,|\, 1\leq i\leq k, 1\leq j\leq d_i-1\}
$$
and the $g$-collection
$$
\{D^{g}_{i,j}\,|\, 1\leq i\leq m, l_i\leq j \leq 2l_i-1,
(i,j)\neq (m,2l_m-1)\}.
$$
Since the linear system $|H|$ is by construction free
(and defines precisely the cover
$\varphi_{|H|}=\sigma \colon V\to{\mathbb P}$), we get by
Corollary 3.2 (inequality (\ref{b4}))
\begin{equation}
\label{b14}
\lambda_{k+1}(x)\leq
\left(
\prod^k_{i=1}\prod^{d_i-1}_{j=2}
\frac{\displaystyle
j+1}{\displaystyle j}
\right)^{-1}\cdot
\left(
\prod^m_{i=1}\prod^{2l_i-1}_{j=l_i}
\frac{\displaystyle
j+1}{\displaystyle j}
\right)^{-1}\cdot
\frac{\displaystyle
2l_m}{\displaystyle 2l_m-1}\cdot
\frac{\displaystyle
a+1}{\displaystyle a},
\end{equation}
where $a=2$, if $\mathop{\rm max}\limits_i \{d_i\}\geq 3$,
and  $a=\mathop{\rm min}\limits_i \{l_i\}$ in the opposite
case. Indeed, to apply the operation $\mathop{\rm min}$ in
the inequality (\ref{b4}) means in the notations of 
Corollary 3.2 to delete from the product
$$
\prod^N_{i=1}\frac{\displaystyle
\mu_i}{\displaystyle a_i}
$$
precisely the $e$ highest factors. In our case these
factors obviously are $k$ twos (corresponding to the
tangent hyperplane sections $D^{f}_{i,1}$) and the
next factor $(a+1)/a$. The factor $2l_m/(2l_m-1)$ comes
into (\ref{b14}) simply  because the divisor 
$D^{g}_{m,2l_m-1}$ is not present in our collection.
Making in (\ref{b14}) the obvious cancellations, we obtain
$$
\lambda_{k+1}(x)\leq
\frac{\displaystyle
2^k}{\displaystyle \mathop{\rm deg}\nolimits V}\cdot
\frac{\displaystyle
2l_m}{\displaystyle 2l_m-1}\cdot
\frac{\displaystyle
3}{\displaystyle 2},
$$
since in any case $a\geq 2$. Finally for the correct
subvariety $Y$ we get:
$$
\frac{\displaystyle \mathop{\rm mult}\nolimits_x}{\displaystyle \mathop{\rm deg}} Y\leq 2^{1-k}\cdot \frac{\displaystyle \mathop{\rm mult}\nolimits_x}{\displaystyle \mathop{\rm deg}} R\leq
$$
$$
\leq 
\frac{\displaystyle
4}{\displaystyle \mathop{\rm deg}\nolimits V}\cdot
\frac{\displaystyle
3l_m}{\displaystyle 4l_m-2}\leq
\frac{\displaystyle
4}{\displaystyle \mathop{\rm deg}\nolimits V},
$$
which is what we need. Proof of the crucial
inequality (\ref{c1}) for a point $x\in V$ such that the point
$p=\sigma (x)\in Q$ is of class 0 is complete.


\subsection{The points of class $e\geq 1$}

Now assume that for the point $x\in V$ the point
$p=\sigma (x)\in Q$ is of class $e\geq 1$, that is,
$$
p=\sigma(x)\in Q\cap W_1\cap\dots\cap W_e.
$$
By the regularity condition the set-theoretic
intersection
$$
\left(
\mathop{\bigcap}\limits_{i,j}D^{f}_{i,j}
\right)\cap
\left(
\mathop{\bigcap}\limits_{i}D^{}_{i}
\right)\cap
\left(
\mathop{\bigcap}\limits_{i,j}D^{g}_{i,j}
\right)
$$
is in a neighborhood of the point $x$ of the
correct codimension
$$
\sum^k_{i=1}(d_i-1)+e+\sum^m_{i=e+1}l_i\geq 1
$$
(the codimension is taken with respect to $V$). Now
let us apply Corollary 3.2 and obtain an upper bound for
$\lambda_2(x)$. Taking into consideration that
$$
\mathop{\rm min}\limits_{
\begin{array}{c}
\scriptstyle
{\cal L}\subset\{1,\dots,N\} \\
\scriptstyle
\sharp {\cal L}=N-c
\end{array}
}\prod_{i\in{\cal L}}\beta_i\times
\mathop{\rm max}\limits_{
\begin{array}{c}
\scriptstyle
{\cal L}\subset\{1,\dots,N\} \\
\scriptstyle
\sharp {\cal L}=c
\end{array}
}\prod_{i\in{\cal L}}\beta_i=\prod^N_{i=1}\beta_i,
$$
we get
$$
\begin{array}{ccc}
\displaystyle
\lambda_2(x)  &  \leq  & \displaystyle
\left(
\prod^k_{i=1}\prod^{d_i-1}_{j=1}
\frac{\displaystyle
j+1}{\displaystyle j}
\right)^{-1}\cdot 2^{-e}\cdot
\left(
\prod^m_{i=e+1}\prod^{2l_i-1}_{j=l_i}
\frac{\displaystyle
j+1}{\displaystyle j}
\right)^{-1}\cdot 4  \\  \\  \displaystyle
   &   &  \|  \\   \\
\displaystyle
  &   &  \displaystyle
\frac{\displaystyle
4}{\displaystyle \mathop{\rm deg}\nolimits V},
\end{array}
$$
which is what we need. Proof of the crucial estimate (\ref{c1})
(and therefore of our theorem) is complete.


{\small

\section*{References}

\noindent
[C1] Corti A., Factoring birational maps of threefolds after Sarkisov. J. Algebraic Geom. {\bf 4} (1995), no. 2, 
223-254.
\vspace{0.5cm} \par \noindent
[C2] Corti A., Singularities of linear systems and 3-fold
birational geometry, in ``Explicit Birational Geometry
of Threefolds'', London Mathematical Society Lecture Note
Series {\bf 281} (2000), Cambridge University Press, 259-312.
\vspace{0.5cm} \par \noindent
[CM] Corti A. and Mella M., Birational geometry of terminal
quartic 3-folds. I, preprint, math.AG/0102096
\vspace{0.5cm}\par\noindent
[CPR] Corti A., Pukhlikov A. and Reid M., Fano 3-fold
hypersurfaces, in ``Explicit Birational Geometry
of Threefolds'', London Mathematical Society Lecture Note
Series {\bf 281} (2000), Cambridge University Press, 175-258.
\vspace{0.5cm} \par \noindent
[CR] Corti A. and Reid M., Foreword to
``Explicit Birational Geometry
of Threefolds'', London Mathematical Society Lecture Note
Series {\bf 281} (2000), Cambridge University Press, 1-20.
\vspace{0.5cm} \par \noindent
[Ch] Cheltsov I. A., Log canonical thresholds on hypersurfaces. Mat. Sb. {\bf 192} (2001), no.8, 155-172 (English translation in Sb. Math. {\bf 192} (2001), no. 7-8, 1241-1257).
\vspace{0.5cm} \par \noindent
[ChPk] Cheltsov I. A. and Park J.,
Global log-canonical thresholds and generalized Eckardt points, Mat. Sb. {\bf 193} (2002), no. 5, 149--160.
\vspace{0.5cm} \par \noindent
[EV] Esnault H. and Viehweg E., Lectures on vanishing
theorems, DMV-Seminar. Bd. {\bf 20.} Birkh\" auser, 1992. \vspace{0.5cm}\par\noindent
[F1] Fano G., Sopra alcune varieta algebriche a tre dimensioni aventi tutti i generi nulli, Atti Acc. Torino {\bf 43} (1908), 973-977.
\vspace{0.5cm} \par \noindent
[F2] Fano G., Osservazioni sopra alcune varieta non razionali aventi tutti i generi nulli, Atti Acc. Torino 
{\bf 50} (1915), 1067-1072.
\vspace{0.5cm} \par \noindent
[F3] Fano G., Nouve ricerche sulle varieta algebriche a tre dimensioni a curve-sezioni canoniche, Comm. Rent. Ac. Sci. {\bf 11} (1947), 635-720.
\vspace{0.5cm} \par \noindent
[Ful] Fulton W., Intersection Theory, Springer-Verlag, 1984.
\vspace{0.5cm} \par \noindent
[G1] Grinenko M. M., Birational automorphisms of a three-dimensional double cone, Mat. Sb. {\bf 189} (1998), 
no. 7, 37-52 (English translation in Sb. Math. {\bf 189} (1998), no. 7-8, 991-1007).
\vspace{0.5cm} \par \noindent
[G2] Grinenko M. M., Birational properties of pencils of del Pezzo surfaces of degrees 1 and 2. Mat. Sb. {\bf 191} (2000), no. 5, 17-38 (English translation in Sbornik: Math. {\bf 191} (2000), no. 5-6, 633-653).
\vspace{0.5cm} \par \noindent
[I1] Iskovskikh V.A., Birational automorphisms of
three-dimensional algebraic varieties, J. Soviet Math.
{\bf 13} (1980), 815-868.
\vspace{0.5cm} \par \noindent
[I2] Iskovskikh, V. A. A simple proof of the non-rationality of a three-dimensional quartic. Mat. Zametki {\bf 65} (1999), no. 5, 667-673 (English translation in Math. Notes {\bf 65} (1999), no. 5-6, 560-564).
\vspace{0.5cm} \par \noindent
[I3] Iskovskikh V.A., Birational rigidity of Fano
hypersurfaces from the viewpoint of Mori theory, Russian Math. Surveys {\bf 56} (2001), no. 2, 3-86.\vspace{0.5cm}
\par\noindent
[IM] Iskovskikh V.A. and Manin Yu.I., Three-dimensional
quartics and counterexamples to the L\" uroth problem,
Math. USSR Sb. {\bf 86} (1971), no. 1, 140-166.
\vspace{0.5cm} \par \noindent
[IP] Iskovskikh V.A. and Pukhlikov A.V., Birational
automorphisms of multi-dimensional algebraic varieties,
J. Math. Sci. {\bf 82} (1996), 3528-3613.
\vspace{0.5cm} \par \noindent
[K] Koll{\'a}r J., et al., Flips and Abundance for
Algebraic Threefolds, Asterisque 211, 1993.
\vspace{0.5cm} \par \noindent
[Kw] Kawamata Y., A generalization of
Kodaira-Ramanujam's vanishing theorem, Math. Ann. {\bf 261}
(1982), 43-46.\vspace{0.5cm}\par\noindent
[M1] Manin Yu. I., Rational surfaces over perfect fields.  Publ. Math. IHES {\bf 30} (1966), 55-113.
\vspace{0.5cm} \par \noindent
[M2] Manin Yu. I. Rational surfaces over perfect fields. II.  Mat. Sb. {\bf 72} (1967), 161-192.
\vspace{0.5cm} \par \noindent
[M3] Manin Yu. I., Cubic forms. Algebra, geometry, arithmetic. Second edition. North-Holland Mathematical Library, {\bf 4.} North-Holland Publishing Co., Amsterdam, 1986.
\vspace{0.5cm} \par \noindent 
[N] Noether M., {\" U}ber Fl{\" a}chen welche Schaaren rationaler Curven besitzen, Math. Ann. {\bf 3} (1871), 161-227.
\vspace{0.5cm} \par \noindent
[P1] Pukhlikov A.V., Birational isomorphisms of four-dimensional quintics, Invent. Math. {\bf 87} (1987), 303-329.
\vspace{0.5cm} \par \noindent
[P2] Pukhlikov A.V., Birational automorphisms of a double space and a double quadric, Math. USSR Izv. {\bf 32} (1989), 233-243.
\vspace{0.5cm} \par \noindent
[P3] Pukhlikov A.V., Birational automorphisms of a three-dimensional quartic with an elementary singularity, Math. USSR Sb. {\bf 63} (1989), 457-482.
\vspace{0.5cm} \par \noindent
[P4] Pukhlikov A.V., A note on the theorem of V.A.Iskovskikh and Yu.I.Manin on the three-dimensional quartic, Proc. Steklov Math. Inst. {\bf 208} (1995), 244-254.
\vspace{0.5cm} \par \noindent
[P5] Pukhlikov A.V., Essentials of the method of maximal
singularities, in ``Explicit Birational Geometry
of Threefolds'', London Mathematical Society Lecture Note
Series {\bf 281} (2000), Cambridge University Press, 73-100.
\vspace{0.5cm} \par \noindent
[P6] Pukhlikov A.V., Birational automorphisms of
three-dimensional algebraic varieties with a pencil of
del Pezzo surfaces, Izvestiya: Mathematics {\bf 62}:1 (1998), 115-155.
\vspace{0.5cm} \par \noindent
[P7] Pukhlikov A.V., Birational automorphisms of Fano
hypersurfaces, Invent. Math. {\bf 134} (1998), no. 2,
401-426. \vspace{0.5cm} \par \noindent
[P8] Pukhlikov A.V., Birationally rigid Fano double
hypersurfaces, Sbornik: Mathematics {\bf 191} (2000), No. 6,
101-126.
\vspace{0.5cm} \par \noindent
[P9] Pukhlikov A.V., Birationally rigid Fano fibrations,
Izvestiya: Mathematics {\bf 64} (2000), 131-150.
\vspace{0.5cm}\par \noindent
[P10] Pukhlikov A.V., Birationally rigid Fano complete
intersections, Crelle J. f\" ur die reine und angew. Math. {\bf 541} (2001), 55-79. \vspace{0.5cm}\par\noindent 
[P11] Pukhlikov A.V., Birationally rigid Fano
hypersurfaces with isolated singularities, Sbornik: Mathematics  {\bf 193} (2002), no. 3,. \vspace{0.5cm}
\par \noindent
[P12] Pukhlikov A.V., Birationally rigid Fano hypersurfaces,
Izvestiya: Mathematics {\bf 66} (2002), no. 6; 
math.AG/0201302.
\vspace{0.5cm}\par \noindent
[R] Reid M., Birational geometry of 3-folds according to
Sarkisov. Warwick Preprint, 1991.
\vspace{0.5cm} \par \noindent
[S1] Sarkisov V.G., Birational automorphisms of conic
bundles, Izv. Akad. Nauk SSSR, Ser. Mat. {\bf 44} (1980),
no. 4, 918-945 (English translation:
Math. USSR Izv. {\bf 17} (1981), 177-202).
\vspace{0.5cm} \par \noindent
[S2] Sarkisov V.G., On conic bundle structures, Izv. Akad.
Nauk SSSR, Ser. Mat. {\bf 46} (1982), no. 2, 371-408
(English translation: Math. USSR Izv. {\bf 20} (1982), no. 2, 354-390).
\vspace{0.5cm} \par \noindent
[S3] Sarkisov V.G., Birational maps of standard ${\mathbb Q}$-Fano fibrations, Preprint, Kurchatov Institute of
Atomic Energy, 1989.
\vspace{0.5cm} \par \noindent
[Sh] Shokurov V.V., 3-fold log flips, Izvestiya:
Mathematics {\bf 40} (1993), 93-202. \vspace{0.5cm}
\par \noindent 
[Sob1] Sobolev I. V., On a series of birationally rigid varieties with a pencil of Fano hypersurfaces. Mat. Sb. 
{\bf 192} (2001), no. 10, 123-130 (English translation in Sbornik: Math. {\bf 192} (2001), no. 9-10, 1543-1551). \vspace{0.5cm} \par \noindent
[Sob2] Sobolev I. V., Birational automorphisms of a class of varieties fibered into cubic surfaces. Izv. Ross. Akad. 
Nauk Ser. Mat. {\bf 66} (2002), no. 1, 203-224.
\vspace{0.5cm} \par \noindent
[V] Viehweg E., Vanishing theorems, Crelle J. f\" ur die reine und angew. Math. {\bf 335} (1982), 1-8.

}

\end{document}